\documentclass[11pt]{article}
\usepackage{amsmath}
\usepackage{amsthm}
\usepackage{graphicx}
\usepackage{amssymb}
\usepackage{amscd}
\theoremstyle{plain}  
  \newtheorem{thm}{Theorem}[section]
  
  \newtheorem{lem}[thm]{Lemma}

\theoremstyle{definition}  
 \newtheorem{defn}[thm]{Definition}
\theoremstyle{remark}
  \newtheorem*{rem}{Remark}

  \newtheorem*{ack}{Acknowledgements}

\newcommand{\C}{\mathbb C} 
\newcommand{\Z}{\mathbb Z} 
  
\newcommand{\R}{\mathbb R}  
\newcommand{\Q}{\mathbb Q}  
\newcommand{\HH}{{\bf H}}
\newcommand{\bHH}{\ov{\bf H}}
\newcommand{\HHD}{(\bHH \times \bHH) \setminus \Delta}
  
\newcommand{\hR}{\hat \R}  
\newcommand{\hQ}{\hat \Q} 
\newcommand{\hC}{\hat \C} 
 
\newcommand{\RR}{\mathcal R} 
\newcommand{\D}{\mathcal D(S)}  
\newcommand{\M}{\mathcal M} 
\newcommand{\B}{\mathcal B}  
\newcommand{\T}{\mathcal T}  
\newcommand{\brho}{\boldsymbol{\rho}}   
\newcommand{\bump}{\mathrm{bump}}  
\newcommand{\exotic}{\mathrm{exotic}}  

\newcommand{\lam}{\lambda}

\newcommand{\fty}{\infty} 
\newcommand{\ep}{\epsilon}

\newcommand{\bd}{\partial} 
\newcommand{\wh}{\widehat}

\newcommand{\ov}{\overline}
\newcommand{\ten}{\cdot}
\newcommand{\sm}{\setminus}
\newcommand{\la}{\langle}
\newcommand{\ra}{\rangle}

\newcommand{\Mod}{\mathrm{Mod}}    

\newcommand{\inte}{\mathrm{int}}

\newcommand{\id}{\mathit{id}}
  
\newcommand{\im}{\mathrm{Im}\,} 
\newcommand{\re}{\mathrm{Re}\,} 

\newcommand{\psl}{{\mathrm{PSL}}_2(\mathbb C)}

\title{Convergence and divergence of \\ 
Kleinian punctured torus groups} 
\author{Kentaro Ito}
\begin{document}

\maketitle

\begin{abstract}
In this paper we give a necessary and sufficient condition 
in which a sequence of Kleinian punctured torus groups converges.   
This result tells us that 
every {\it exotically} convergent sequence 
of Kleinian punctured torus groups 
is obtained by the method due to Anderson and Canary 
(Invent.~Math.~1996). 
Thus we obtain a complete description of  
the set of points at which the space 
of Kleinian punctured torus groups self-bumps. 
We also  discuss Hausdorff limits of 
sequences of Bers slices. 
\end{abstract}

\section{Introduction}
One of the central issues 
in the theory of Kleinian groups is to understand 
the structures of deformation spaces of Kleinian groups. 
In this paper we consider Kleinian punctured torus groups, 
one of the simplest classes of 
Kleinian groups with a non-trivial deformation theory. 
Let $S$ be a once-punctured torus. 
The deformation space ${\cal D}(S)$  
of Kleinian punctured torus groups  
is the space of conjugacy classes $[\rho]$ of discrete faithful representations 
$\rho:\pi_1(S) \to \psl$ 
which takes a loop surrounding the cusp to a parabolic element. 
Although the interior of $\D$ 
is parameterized by a product of Teichm\"{u}ller spaces of $S$, 
its boundary is quite complicated. 
For example, McMullen \cite{Mc1} showed that $\D$ self-bumps by using 
the method developed by Anderson and Canary \cite{AC} 
(see also \cite{BH}): 
here we say that $\D$ {\it self-bumps}
if there is a point $[\rho]$ on the boundary 
such that for any sufficiently small neighborhood of $[\rho]$, 
the intersection of the neighborhood 
with the interior of $\D$ is disconnected. 
Furthermore,  Bromberg \cite{Brom} recently showed 
that $\D$ is not even locally connected. 
We refer the reader to \cite{Ca1} and \cite{Ca2} 
for more information on the topology of 
deformation spaces of general Kleinian groups. 

In this paper we characterize 
sequences in $\D$ which give rise to the self-bumping of $\D$. 
Roughly speaking, all such a sequence is obtained by the 
construction developed by Anderson and Canary \cite{AC}. 
To describe our results, we review the basic setting. 

Let $\T(S)$ denote Teichm\"{u}ller space 
of once punctured torus $S$. 
We denote by $\ov{\T}(S)$ 
the Thurston compactification of $\T(S)$ 
with the set $\mathcal{PL}(S)$ of projective measured laminations on $S$.  
Then Minsky's ending lamination theorem for punctured torus groups \cite{Mi} assert that 
all elements $[\rho] \in \D$ are classified by their end invariants 
$\nu([\rho])  \in \left(\ov{\T}(S) \times \ov{\T}(S)\right)  \sm \Delta$,  
where $\Delta$ is the diagonal of $\mathcal{PL}(S) \times \mathcal{PL}(S)$.  
Moreover, he showed that the inverse 
\begin{eqnarray*}
Q=\nu^{-1}: \left(\ov{\T}(S) \times \ov{\T}(S)\right)  \sm \Delta \to \D
\end{eqnarray*}
of the map $\nu$ is bijective and continuous.

The purpose of this paper is to obtain a necessary  and 
sufficient condition of sequences $\{(x_n,y_n)\}$ in 
$\left(\ov{\T}(S) \times \ov{\T}(S)\right) \sm \Delta$ 
 such that $\{Q(x_n,y_n)\}$ converges in $\D$. 
Since the map $Q$ is continuous on 
$\left(\ov{\T}(S) \times \ov{\T}(S)\right) \sm \Delta$,  
we concentrate our attention to  the behavior of the sequence 
$\{Q(x_n,y_n)\}$ such that $\{(x_n,y_n)\}$ converges to 
a point $(x_\fty,x_\fty)$ in $\Delta$. 
If such a sequence $\{Q(x_n,y_n)\}$  converges, we say that it is an  
{\it exotically} convergent sequence.  

It was shown by Ohshika \cite{Oh1} that if 
$x_\fty \in \mathcal{PL}(S)$ is not a simple closed curve,  then there is no 
exotically convergent sequence in $\D$ 
both of whose end invariants converge to $x_\fty$. 
On the other hand, in the case that  $x_\fty$ is a simple closed curve,  
there is an exotically convergent sequence 
both of whose end invariants converge to $x_\fty$. 
Such a sequence was first obtained by McMullen \cite{Mc1} 
by using the method due to Anderson and Canary \cite{AC}: 
Let $c$ be a simple closed curve on $S$ and 
let $\tau$ denote the Dehn twist around $c$. 
Then, for given $x,y \in \T(S)$ and an integer $p$ with $p \ne 0,-1$, 
the sequence 
\begin{eqnarray*}
\{Q(\tau^{p n} x, \tau^{(p+1)n} y)\}_{n=1}^\fty 
\end{eqnarray*}
converges in $\D$, whereas 
the sequence 
$$
\{(\tau^{p n} x, \tau^{(p+1)n} y)\}_{n=1}^\fty
$$
 in $\left(\ov{\T}(S) \times \ov{\T}(S)\right) \sm \Delta$ 
converges to $(c,c) \in \Delta$.  

Our main result,  Theorem \ref{main} below 
(see Theorems \ref{horo} and \ref{tan}),  implies that 
all  exotically convergent sequences 
are essentially obtained by the method of Anderson and Canary.  
More precisely, it states that  for a given simple closed curve $c$ on $S$,  
if either $\{x_n\}$ or $\{y_n\}$ converge ``horocyclically"  
to $c$  in $\ov \T(S)$  then $\{Q(x_n,y_n)\}$ diverges, 
and that if both $\{x_n\}$ and $\{y_n\}$ converge 
``tangentially" to $c$ then $\{Q(x_n,y_n)\}$ converges provided that 
the speeds of their convergence to $c$ are in the ratio of 
$p:p+1$ for some integer $p$. 

We denote by $l_x(c)$ the hyperbolic length of  the geodesic realization 
of a simple closed curve $c$ on the Riemann surface $x \in \T(S)$. 

\begin{thm}[Theorems \ref{horo} and \ref{tan}]\label{main}
Suppose that a sequence 
$\{(x_n,y_n)\}$ in $(\ov\T(S) \times \ov\T (S)) \sm \Delta$ 
converges to $(c,c) \in \Delta$ for some simple closed curve $c$ on $S$
as $n \to \fty$. 
Then we have the following: 
\begin{enumerate}
    \item If either $\{l_{x_n}(c)\}$ or  
    $\{l_{y_n}(c)\}$ tend to zero,  the sequence $\{Q(x_n, y_n)\}$ diverges in $\D$.  
  \item Suppose that both 
 $\{l_{x_n}(c)\}$ and 
    $\{l_{y_n}(c)\}$  are uniformly bounded below by a positive constant.    
We may also assume, by pass to a subsequence if necessary, that 
there exist sequences $\{k_n\}$, $\{l_n\}$ of 
    integers such that both $\{\tau^{k_n}x_n\}$, $\{\tau^{l_n}y_n\}$  
    converge 
    in $\ov\T(S) \sm \{c\}$, where $\tau$ is the Dehn twist around $c$. 
   In this situation,  the sequence $\{Q(x_n,y_n)\}$ converges in $\D$ 
    if and only if there exist an integer $p$ 
    such that 
\begin{eqnarray*}
    (p+1)k_n-pl_n \equiv \text{const.}   
\end{eqnarray*}
    for all $n$ large enough. 
\end{enumerate}
\end{thm}

\begin{rem}
Recently, Ohshika \cite{Oh2} obtained a generalization of Theorem \ref{main} 
to general hyperbolic surfaces of finite area. 
\end{rem}

One of our main tool in the proof of Theorem \ref{main} are the ``re-marking trick" on representations 
which is originally due to Kerckhoff and Thurston \cite{KT} (see also \cite{Brock}). 
We also make essential use of Minsky's Pivot Theorem \cite{Mi} 
(see Theorem \ref{pivot}) which 
relates a pair of end invariants of an element in $\D$ to 
the complex translation length of a 
short closed geodesic in the quotient manifold.

We now give a brief sketch of  the proof of Theorem \ref{main} (2). 
For simplicity we assume that $\tau^{k_n}x_n$ and $\tau^{l_n}y_n$ are constantly equal 
to $u,v \in \ov{\T}(S) \sm \{c\}$, respectively. 
Then we have 
\begin{equation*}
Q(x_n,y_n)=Q(\tau^{-k_n} u,\tau^{-l_n}v)
=Q(u,\tau^{k_n-l_n}v) \circ \tau_*^{k_n}, \tag{1.1}
\end{equation*}
where $\tau_*:\pi_1(S) \to \pi_1(S)$ is the isomorphism 
induced by the Dehn twist $\tau$.   
As we will see in Section 4, we only need to consider the case 
that $|l_n-k_n|$ diverge as $n \to \fty$.  
In this case we have 
$$
Q(u,\tau^{l_n-k_n}v) \to Q(u,c)
$$ 
in $\D$.  
Let $\rho_n:\pi_1(S) \to \psl$ be representatives in the conjugacy classes  
$Q(u,\tau^{l_n-k_n}v)$ such that $\{\rho_n\}$ converges algebraically.  
Then $\{\rho_n(c)\}$ converges to a parabolic transformation $\delta$. 
Moreover, we can  show that 
the sequence $\{\la\rho_n(c)\ra\}$ of the cyclic groups  converges geometrically to 
a rank-2 parabolic group $\la \delta, \hat \delta \ra$ 
with  $\hat \delta=\lim_{n \to \fty}\rho_n(c)^{l_n-k_n}$ (see Theorem \ref{prime}). 

On the other hand,  
since $\{Q(u,\tau^{k_n-l_n}v) \}$ converges, 
one can see from (1.1) that 
$\{Q(x_n,y_n)\}$ converges if and only if 
$\{\rho_n(c)^{k_n}\}$ converges (see Lemma \ref{trick}). 
If $\{\rho_n(c)^{k_n}\}$ converges,  its limit must lie in the rank-2 parabolic group 
 $ \la \delta, \hat \delta\ra$. 
This implies that there exists integers $p,\,q$ such that 
$k_n \equiv p(l_n-k_n)+q$, or,  equivalently  
$(p+1)k_n-pl_n \equiv q$.  
The converse is easier.

This paper is organized as follows. 
In Section 2, we give the basic notion and definitions.  
Especially we recall Minsky's ending lamination theorem for punctured torus groups. 
In Section 3, we study the relation between  
geometric limits of loxodromic cyclic groups 
and the complex translation lengths of generators.   
Using preparations of Section 3, 
we prove the main theorem in Section 4. 
Applications of the main theorem are discussed in 
Sections 5, 6 and 7: 
In Section 5, we recall Bers and Maskit slices, and the Maskit embedding. 
We rephrase the main theorem in terms of the Maskit embedding. 
In Section 6, we consider Hausdorff limits of sequences of Bers slices 
and obtain a condition in which the Hausdorff limit 
of a sequence of Bers slices 
is strictly bigger than 
the limit of point wise convergence (see Theorem \ref{bersgeom}). 
In Section 7, we obtain a complete description of  
the set of points at which $\D$ self-bumps (see Theorem \ref{bump}). 

\begin{ack}
The author would like to thank 
Ken'ichi Ohshika and Hideki Miyachi for useful conversations and information. 
The author also thank the referee's valuable comments and suggestions. 
\end{ack}

\section{Preliminaries}

\subsection{Teichm\"{u}ller space}

Let $S$ be a once-punctured torus.  
Throughout this paper, 
we fix an ordered pair $(\alpha, \beta)$ of generators for $\pi_1(S)$, 
which determines a positively oriented ordered pair 
$([\alpha], [\beta])$ of generators for $H_1(S)$. 
Here $[\alpha]$ and $[\beta]$ denote 
homology classes of $\alpha$ and $\beta$, respectively. 

The Theichm\"{u}ller space ${\cal T}(S)$ consists of pairs $(f,X)$, 
where $X$ is a hyperbolic Riemann surface of finite area and 
$f:S \to X$ is an orientation preserving homeomorphism.  
Two pairs $(f_1,X_1)$ and $(f_2,X_2)$ represent the same point in ${\cal T}(S)$ 
if there is a holomorphic isomorphism $g:X_1 \to X_2$ 
such that $g \circ f_1$ is isotopic to $f_2$. 

The Teichm\"{u}ller space 
${\cal T}(S)$ is identified with the upper half plane 
$\HH=\{z \in \C \,|\, \mathrm{Im}\,z >0\}$ as follows: 
To a point $x$ in $\HH$, we associate the point $(f,X) \in {\cal T}(S)$ where 
$X$ is the quotient $\C/(\Z \oplus \Z x)$ 
with one point removed, and  
 $f:S \to X$ is an orientation preserving homeomorphism
which take the curves $[\alpha]$ and $[\beta]$ 
to the images of the segments $[0,1]$ and $[0,x]$ in $X$, respectively.  

In this identification, 
Thurston's compactification $\ov{\cal T}(S)$ of the Teichm\"{u}ller space 
${\cal T}(S)$ with the set ${\cal PL}(S)$ of projective measured laminations on $S$ 
corresponds to the closure 
$\ov \HH=\HH \cup \hR$ of $\HH$ in the Riemann sphere $\hC$, 
where $\hR=\R \cup \{\fty\}$. 
Then the simple closed curve on $S$ that represent an unoriented homology class 
$\pm(s[\alpha]+t[\beta]) \in H_1(S)$ is identified with the rational number 
$-s/t \in \hQ=\Q \cup \{\fty\}$. 
Especially, the homology classes 
$[\alpha], \,[\beta]$ and $[\alpha^{-1}\beta]$  
correspond to $1/0=\fty,\, 0$ and $1$ in $\hQ$, respectively. 

We let $l_x(c)$ denote the hyperbolic length 
of the geodesic in the homotopy class 
of a simple closed curve $c$ on a Riemann surface $x \in \T(S)$.  
Via the identification $\ov \T(S)$ with $\bHH$ described above, 
a horocircle in $\HH$ touching $\bd \HH$ at $c \in \hQ$ 
is the set of Riemann surfaces 
$x \in \HH=\T(S)$  
with the same hyperbolic lengths $l_x(c)$. 

\subsection{Kleinian groups}

A {\it Kleinian group} $\Gamma$ is a discrete subgroup of $\psl$, 
which acts on the hyperbolic 3-space $\HH^3$ as isometries, and on 
the sphere at infinity $\bd \HH^3=\hat \C$ 
as conformal automorphisms. 
The {\it region of discontinuity} 
$\Omega_\Gamma$ for a Kleinian group $\Gamma$ is the largest open subset of 
$\wh{\C}$ on which $\Gamma$ acts properly discontinuously, 
and the {\it limit set} $\Lambda_\Gamma$ of $\Gamma$ 
is its complement $\wh{\C}-\Omega_\Gamma$.  

We say that a sequence $\{\Gamma_n\}_{n=1}^\fty$ 
of Kleinian groups converges {\it geometrically} to a subgroup 
$\hat \Gamma$ of $\psl$ if 
$\Gamma_n \to \hat \Gamma$ in the Hausdorff topology on closed subsets of 
$\psl$,  
or more precisely, if the following conditions are satisfied: 
\begin{enumerate}
\item for any $\gamma \in \hat\Gamma$ there exists a sequence $\gamma_n \in \Gamma_n$  
such that $\gamma_n \to \gamma$, and 
\item if elements $\gamma_{n_j} \in \Gamma_{n_j}$ converge to $\gamma$, 
then $\gamma \in \hat \Gamma$. 
\end{enumerate}

Let $G$ be an abstract group and let $\rho_n:G \to \Gamma_n \subset \psl$ 
be a sequence of  homomorphisms 
form $G$ onto Kleinian groups $\Gamma_n$. 
We say that $\{\rho_n\}$ converges  {\it algebraically} to a 
homomorphism $\rho_\fty:G \to \Gamma_\fty \subset \psl$ 
if $\rho_n(g) \to \rho_\fty(g)$  in $\psl$ for every $g \in G$. 
In this situation, 
 if $\{\Gamma_n\}$ 
also converges geometrically to $\Gamma_\fty$, 
we say that $\{\Gamma_n\}$ converges {\it strongly} 
to $\Gamma_\fty$. 
 
Let $\RR(S)$ denote the space of conjugacy classes $[\rho]$ of
representations $\rho$  of $\pi_1(S)$ into $\psl$ which 
take the commutator of $\alpha$ and $\beta$ to parabolic elements.  
We endow $\RR(S)$ the algebraic topology; i.e., 
a sequence $[\rho_n] \in \RR(S)$ converges to  
$[\rho] \in \RR(S)$ if and only if there exist representatives 
$\rho_n \in [\rho_n]$ and $\rho \in [\rho]$ 
such that $\rho_n \to \rho$ algebraically.  
It is known that $\RR(S)$ has a structure of 
2-dimensional complex  manifold,  and 
especially it is a Hausdorff space.  
More precisely, $\RR(S)$ can be identified with the space 
$$
\{(x,y,z) \in \C^3 \sm \{(0,0,0)\} : x^2+y^2+z^2=xyz\}/\Z_2 \times \Z_2, 
$$
where $\Z_2 \times \Z_2$ acts by changing signs 
of two of three entries (see \cite{Bow}). 

Let $\D \subset \RR(S)$ denote the set of 
discrete faithful representations. 
It is known that $\D$ is a closed subset of $\RR(S)$ (see \cite{Jo}).  

\subsection{Minsky's ending lamination theorem} 

Here we review Minsky's work (\cite{Mi}) on the classification of 
elements $[\rho] \in \D$ by using pairs of  their {\it end invariants}  
$(\nu_-,\nu_+) \in \HHD$, where 
$\Delta$ is the diagonal of 
$\bd \HH \times \bd \HH$: 
\begin{eqnarray*}
\Delta=\{(x,y) \in \bd \HH \times \bd \HH \,|\, x=y \}. 
\end{eqnarray*}
We refer the reader to \cite{Mi} for more details. 

Given $[\rho] \in \D$ with $\Gamma=\rho(\pi_1(S))$, 
Bonahon's Theorem \cite{Bo} guarantees that 
there exists an orientation preserving homeomorphism 
$$
\varphi:S \times (-1,1) \to \HH^3/\Gamma
$$
which induces the representation $\rho$. 
The region of discontinuity 
$\Omega_\Gamma$ (possibly empty) decomposes 
into two parts $\Omega_\Gamma^+$, $\Omega_\Gamma^-$ 
such that $\Omega_\Gamma^+/\Gamma$ and $\Omega_\Gamma^-/\Gamma$ 
are the limits of $\varphi(S \times \{t\})$ in 
$(\HH^3 \cup \Omega_\Gamma)/\Gamma$ 
as $t \to  \pm 1$, respectively.  Let $s$ denote either $+$ or $-$. 
Then $\Omega_\Gamma^s/\Gamma$ is either once-punctured torus, 
thrice-punctured sphere or empty. 
For each cases, we define the end invariant $\nu_s \in \ov \T(S)=\bHH$ 
for $[\rho]$ as follows: 

\begin{enumerate}
\item If $\Omega_\Gamma^s/\Gamma$ is a once-punctured torus, 
its conformal structure with the marking induced by $\varphi$ 
determines a point $\nu_s$ in $\T(S)=\HH$. 
\item If  $\Omega_\Gamma^s/\Gamma$ is  a thrice-punctured sphere, 
there is a simple closed curve $c$ on $S$ 
such that $\varphi$ extends continuously to 
$(S \sm c) \times \{s \cdot 1\}$ as 
a homeomorphism onto $\Omega_\Gamma^s/\Gamma$. 
In this case we let $\nu_s=c \in \hQ \subset \bd \HH$.  
\item  If $\Omega_\Gamma^s/\Gamma$ is empty, 
there is a sequence  $c_n \in \hQ$ of simple closed curves on $S$ such that  
their geodesic realizations $c_n^*$ in $\HH^3/\Gamma$ diverge to 
the end of $\HH^3/\Gamma$ associated to the sign $s$
and that $c_n \in \hQ$ converges in $\hR$ to some irrational number $c_\fty$. 
In this case we let $\nu_s=c_\fty \in \hR \sm \hQ \subset \bd \HH$.  
\end{enumerate}

Thus we obtain a map $\nu:[\rho] \mapsto (\nu_-,\nu_+)$ from $\D$ to 
$\bHH \times \bHH$. 
It is known that 
the map $\nu$ is surjective onto $\HHD$ 
as a consequence of Bers' Simultaneous Uniformization Theorem, 
and of Thurston's Double Limit Theorem. 
Minsky \cite{Mi} showed that the map $\nu$ is injective; that is, 
all elements $[\rho] \in \D$ are classified by their end invariants.  
Furthermore he showed that the inverse 
$Q:=\nu^{-1}$ of the map $\nu:[\rho] \mapsto (\nu_-,\nu_+)$ is continuous: 
\begin{thm}[Minsky's ending lamination theorem \cite{Mi}]\label{Min}
The map 
$$  
Q: \HHD \to \D
$$ 
is a continuous bijection. 
\end{thm}
This theorem is an extension of the following well-known result: 
\begin{thm}[Ahlfors, Bers, Kra, Maskit, Marden and Sullivan]
The map 
$$
Q:  \HH \times \HH \to \inte(\D)  
$$
is a homeomorphism onto the interior $\inte(\D)$ of $\D$. 
\end{thm}

We remark that the map $\nu:\D \to \HHD$ 
is not continuous at the boundary of $\D$. 
In fact, it was shown by McMullen  \cite{Mc1} by using the method of 
Anderson and Canary \cite{AC} 
that there is a convergent sequence in $\D$ whose end invariants 
converge to a point in $\Delta$.  
Therefore the map $Q$ is not a homeomorphism. 
We will explain more details of this phenomenon in Section 4.

\subsection{Action of $\Mod(S)$} 

In this subsection, we recall the action of the mapping class group 
on the space  $\RR(S)$ of representations.  

The {\it mapping class group} $\Mod(S)$ 
is the group of isotopy classes of 
orientation-preserving homeomorphisms from $S$ to itself. 
The action of $\sigma \in \Mod(S)$  on 
${\cal T}(S)$ 
is defined by 
$$
\sigma(f,X):=(f \circ \sigma^{-1},X). 
$$
Via our identification of $\T(S)$ with $\HH$, 
$\Mod(S)$ is identified with $\mathrm{PSL}_2(\Z)$. 
The action of $\Mod(S)=\mathrm{PSL}_2(\Z)$  on $\T(S)=\HH$ 
naturally extends to automorphisms of 
$\ov{\T}(S)=\bHH$. 

The action of $\sigma \in \Mod(S)$ on a representation 
$\rho:\pi_1(S) \to \psl$ is defined by
$$
\sigma \ten \rho:=\rho \circ \sigma_*^{-1}, 
$$
where $\sigma_*$ is the group automorphism of $\pi_1(S)$   
induced by $\sigma$. 
Then the action of $\sigma  \in \Mod(S)$ on  $\RR(S)$ is defined by 
$$
\sigma \ten [\rho]:=[\sigma \ten \rho]. 
$$
Note that these actions are compatible with the map 
$Q$, i.e., we have 
$$
\sigma \ten Q(x,y)=Q(\sigma x, \sigma y)
$$
for every $\sigma \in \Mod(S)$ and $(x,y) \in \HHD$. 

Throughout of this paper, 
we denote by $\tau \in \Mod(S)$ the Dehn twist around the curve $[\alpha]$, 
whose orientation is chosen so that 
the group isomorphism $\tau_*:\pi_1(S) \to \pi_1(S)$ satisfies 
$$
\tau_*(\alpha)=(\alpha) \quad \text{and} \quad \tau_*(\beta)=\alpha^{-1} \beta. 
$$
Then $\tau$ acts on $\HH={\cal T}(S)$ by 
$$
z \mapsto z+1.  
$$
Given a representation $\rho:\pi_1(S) \to \psl$, 
one can check that 
$$
(\tau \ten \rho) (\alpha)=\rho(\alpha) \quad \text{and}\quad 
(\tau \ten \rho) (\beta)=\rho(\alpha \beta). 
$$
It is often useful to consider the restriction of a representation 
$\rho:\pi_1(S) \to \psl$ to the subgroup  $H:=\pi_1(S \sm [\alpha])$ of $\pi_1(S)$. 
Here $H$ is a rank-$2$ free subgroup of $\pi_1(S)$
generated by $\alpha$ and $\beta^{-1}\alpha\beta$. 
The following lemma is clear from a geometric point of view: 
\begin{lem}\label{H}
$(\tau \ten \rho)|_H =\rho|_H$. 
\end{lem}

\begin{proof}
Since $(\tau \ten \rho)(\alpha)=\rho(\alpha)$, 
it suffices to  show that 
$(\tau \ten \rho)(\beta^{-1} \alpha \beta)=\rho(\beta^{-1} \alpha \beta)$. 
This is directly verified as follows: 
\begin{eqnarray*}
(\tau \ten \rho)(\beta^{-1} \alpha \beta)&=&(\tau \ten \rho)(\beta)^{-1} (\tau \ten \rho)(\alpha) (\tau \ten \rho)(\beta) \\
&=& \rho(\alpha \beta)^{-1} \rho(\alpha) \rho(\alpha\beta) \\
&=& \rho (\beta^{-1} \alpha \beta). 
\end{eqnarray*}
\end{proof}

\section{Complex translation lengths and limits of cyclic groups}

In this section, we study 
geometric limits of loxodromic cyclic groups 
when their generators converge to parabolic transformations. 
The result is applied in the next section to prove our main theorems. 

Recall that the {\it complex translation length}  
$\lambda(\gamma)$ of a loxodromic element $\gamma \in \psl$ 
is defined to be 
$$
\lambda(\gamma)=\log(\gamma'(z))=l+i\theta 
$$ 
with  $l>0$ and $ \theta \in (-\pi,\pi]$,   
where $\gamma'(z)$ is the derivative at the repelling fixed point of $\gamma$.  
This is equivalent to say that $\gamma$ is conjugate in $\psl$ to 
$$
\left(\begin{array}{cc}
e^{\lam(\gamma)/2} & 0 \\ 
0 & e^{-\lam(\gamma)/2} 
\end{array}\right). 
$$
The number $\lam(\gamma)$, or  its reciprocal  $2\pi i/\lam(\gamma)$,  is important to understand the 
modulus of the quotient  torus  $\Omega(\la \gamma \ra)/ \la \gamma \ra$: 
\begin{eqnarray*}
\Omega(\la \gamma \ra)/\la \gamma \ra 
&\cong&(\C \sm \{0\})/\la e^{\lam(\gamma)}z \ra \\
&\cong&\C/\la z+\lam(\gamma),z+2\pi i \ra \\
&\cong&\C  / \la z+1,z+2\pi i/\lam(\gamma)\ra. 
\end{eqnarray*}

For  a given representation $[\rho]=Q(x,y)$ in $\D$, 
if $\rho(\alpha)$ is close to a parabolic transformation, 
the complex translation length of $\rho(\alpha)$ is 
related to the end invariants $x,y \in \bHH$ of $[\rho]$ as follows: 

\begin{thm}[Minsky's Pivot Theorem  \cite{Mi}]\label{pivot} 
There exist positive constants $\ep$, $C$ 
which satisfy the following: 
Given an element $[\rho]=Q(x,y)$ in $\D$, 
if the real part of $\lam(\rho(\alpha))$ is less than $\ep$ 
then 
\begin{eqnarray*}
d_{\HH}\left(\frac{2\pi i}
{\lam(\rho(\alpha))}, x -\ov{y}+i \right)<C, 
\end{eqnarray*}
where $d_{\HH}$ denotes the hyperbolic distance in  $\HH$. 
\end{thm}

Note that if a sequence $\{\gamma_n\}$ in $\psl$ converges  to a parabolic transformation, 
then $\lam(\gamma_n) \to 0$. 
We characterize two ways of convergence 
$\lam(\gamma_n) \to 0$ as follows (see Figure 1): 
\begin{defn}\label{horotan0}
Suppose that a sequence $\{\lam_n\}$ in the right-half plane 
$\C_+=\{z \in \C \,|\, \re z>0\}$ converges to $0$. 
We say that $\lam_n \to 0$ \it{horocyclically} if for any $\epsilon>0$,  
$|\lam_n -\epsilon|<\epsilon$ for all large $n$, and that
 $\lam_n \to 0$ \it{tangentially} if there is a constant $\epsilon_0>0$ 
such that   
$|\lam_n -\epsilon_0|>\epsilon_0$ for all  $n$.  
\end{defn}
\begin{figure}
\begin{center}
\includegraphics[width=8cm]{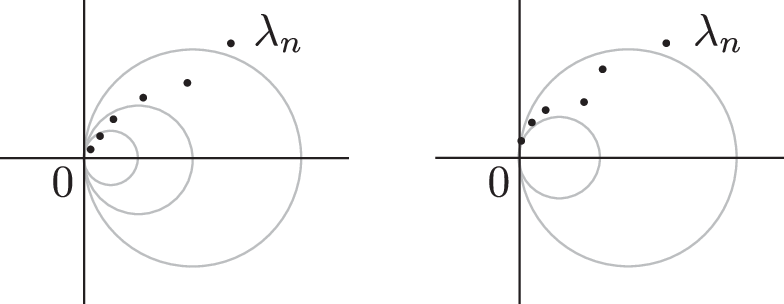}
\caption{Horocyclic convergence (left) and tangential convergence (right). }
\end{center}
\end{figure}
Note that $\lam_n \to 0$ horocyclically if and only if $|\im(2\pi i/\lam_n)| \to \fty$, 
and that tangentially if and only if $|\im(2\pi i/\lam_n)|$ are uniformly bounded above. 

In the proofs of our main theorems, we will make essential use of the following 
theorem, whose essence can be found in McMullen \cite[Theorem 5.1]{Mc2}. 
Our argument is based on that of J{\o}rgensen and Marden \cite[Section 5]{JM} 
and Marden \cite[Section 4.9]{Ma}.

\begin{thm}\label{cyclic}
Suppose that a sequence  $\{\gamma_n\}$  of 
loxodromic transformations converges to a parabolic transformation 
$\delta$ in $\psl$. 
 Then we have the following: 
\begin{enumerate}
\item 
Suppose that $\lam(\gamma_n) \to 0$ horocyclically. 
Then $\la \gamma_n \ra \to \la \delta \ra$ strongly.  

\item  Suppose that $\lam(\gamma_n) \to 0$ tangentially. 
In addition, by pass to a subsequence if necessary, 
we may assume  that there is a sequence 
$\{m_n\}$ of integers such  that 
\begin{equation*}
\frac{2\pi i}{\lam(\gamma_n)}-m_n
\end{equation*}
converges to some $\omega \in \C$. 
In this situation, if  $\omega \not\in \R$,  
the sequence $\la \gamma_n \ra$ converges  geometrically to a 
rank-$2$ parabolic group $\la \delta,  \hat \delta \ra$ 
generated by $\delta=\lim_{n \to \fty} \gamma_n$ and 
$\hat \delta=\lim_{n \to \fty} \gamma_n^{m_n}$. 
Furthermore,  there is an element in $\psl$ which conjugates 
$\delta, \, \hat \delta$ to translations $z+1,\,z-\omega$, respectively. 
\end{enumerate}
\end{thm}

\begin{proof}
We denote $\lam(\gamma_n)$ by $\lam_n$ for simplicity. 
We may assume that $\delta(z)=z+1$. 
Then,  as $n$ tend to $\fty$, $\gamma_n(0)$ tend to $1$ and 
the attracting fixed points $a_n$ of $\gamma_n$ tend to $\fty$.  
Now let $\theta_n \in \psl$ be a M\"{o}bius transformation 
which takes $0,\,\gamma_n(0)$ and $a_n$ to $0,\,1$ and $\fty$, respectively.  
Then $\theta_n \gamma_n \theta_n^{-1}$ is of the form 
$\theta_n \gamma_n \theta_n^{-1} (z)=\alpha_nz+1$.  
By considering the complex translation length of $\gamma_n$, 
one can see  that $\alpha_n=e^{\lam_n}$. 
Since $\theta_n \to \id$ in $\psl$, the geometric limit of the sequence
$\la \gamma_n \ra$ is equal to that of 
$\la \theta_n \gamma_n \theta_n^{-1}\ra$. 
Therefore we may assume that $\gamma_n$ is of the form 
\begin{eqnarray*}
\gamma_n(z)=e^{\lam_n}z+1
\end{eqnarray*}
without loss of generality. 
Then one can easily check  that 
\begin{equation}
\gamma_n^{k_n}(z)=e^{\lam_n k_n} z+\frac{e^{\lam_n k_n}-1}{e^{\lam_n}-1}.  \tag{3.1}
\end{equation}

(1) 
To show that $\la \gamma_n \ra \to \la \delta \ra$ strongly, 
suppose for contradiction that there is a subsequence of $\{\gamma_n\}$, 
which is also denoted by $\{\gamma_n\}$, 
 and a sequence 
$\{k_n\}$ of integers with $|k_n| \to \fty$ 
such that $\{\gamma_{n}^{k_n}\}$ converges in $\psl$. 
Then the sequence 
\begin{equation}
\frac{e^{\lam_n k_n}-1}{e^{\lam_n}-1} \tag{3.2}
\end{equation} 
converges in $\C$. 
Especially $e^{\lam_n k_n} \to 1$, and hence 
there is a sequence of integers $\{p_n\}$ such that 
$\lam_n k_n-2 p_n \pi i \to 0$. 
By Taylor's formula for $e^x$,  the limit of   (3.2) is equal to the limit of 
\begin{equation*}
 \frac{\lam_n k_n-2p_n \pi i}{\lam_n}
=k_n-p_n\frac{2\pi i}{\lam_n}.  
\end{equation*} 
But the sequence $k_n-p_n 2\pi i /\lam_n$ diverges in $\C$ 
from our assumption that $|k_n| \to \fty$ and $|\im(2\pi i/\lam_n)| \to \fty$. 
This is the desired contradiction. 
Therefore we have shown that  $\la \gamma_n\ra \to \la \delta \ra$ strongly. 

(2)  
We first show that $\{\gamma_n^{m_n}\}$ converges in $\psl$ to 
 $\hat \delta(z)=z-\omega$. 
Note from (3.1) that 
\begin{equation*}
\gamma_n^{m_n}(z)=e^{\lam_n m_n} z+\frac{e^{\lam_n m_n}-1}{e^{\lam_n}-1}.  
\end{equation*}
Since $2 \pi i/\lam_n -m_n \to \omega$ and $|m_n| \to \fty$, one can see that 
$2 \pi i/(\lam_n m_n) \to 1$, and thus $\lam_n m_n \to 2\pi i$. 
By Taylor's formula for $e^x$,  we have 
\begin{eqnarray*}
\lim_{n \to \fty}\frac{e^{\lam_n m_n}-1}{e^{\lam_n}-1}
=\lim_{n \to \fty}\frac{\lam_n m_n-2 \pi i}{\lam_n}
=\lim_{n \to \fty}\left(m_n-\frac{2\pi i}{\lam_n}\right)=-\omega. 
\end{eqnarray*}
Therefore 
$\gamma_n^{m_n}(z) \to \hat \delta(z)= z-\omega$ uniformly on any compact subset of $\C$, 
and hence $\gamma_n^{m_n} \to \hat \delta$ in $\psl$. 

We next show that $\gamma_n \to \la \delta,\hat \delta\ra$ strongly. 
 Suppose that  
 there is a subsequence of $\{\gamma_n\}$, 
which is also denoted by $\{\gamma_n\}$, 
 and a sequence 
$\{k_n\}$ of integers with $|k_n| \to \fty$ 
such that $\{\gamma_{n}^{k_n}\}$ converges in $\psl$. 
We have to show that the limit lies in $\la \delta,\hat \delta\ra$.
By the same argument as in (1),  
we see that the sequence 
(3.2) converges to some $\zeta \in \C$, 
and that 
there is a sequence $\{p_n\}$ of integers such that 
 \begin{equation}
\zeta=\lim_{n \to \fty} \left(k_n-p_n\frac{2\pi i}{\lam_n}\right)
=\lim_{n \to \fty} \left(-p_n \left( \frac{2\pi i}{\lam_n}-m_n\right)+k_n-p_n m_n\right).  \tag{3.3}
\end{equation}
Since $2\pi i/\lam_n-m_n \to \omega$ and $\im \omega \ne 0$, 
we obtain $p_n \equiv p$ for all large $n$  by 
considering the imaginary part of  the right-most side of  (3.3). 
We next consider the real part of the right-most side of (3.3) to see that 
 $k_n -p m_n \equiv q$ for some integer $q$. 
 Thus we have $k_n \equiv p m_n +q$ and hence 
$\gamma_n^{k_n} =(\gamma_n^{m_n})^p \gamma_n^q \to \hat \delta^p \delta^q$. 
Therefore the limit of $\gamma_n^{k_n}$ lies in $\la \delta, \hat \delta\ra$. 
\end{proof}

\section{Conditions for convergence and divergence}

In this section, we will prove our main theorems; 
Theorem \ref{horo} and Theorem \ref{tan}.

Suppose that we are given a sequence 
$$
\{(x_n,y_n)\}_{n=1}^\fty
$$ in $\HHD$ 
which converges to some $(x_\fty, y_\fty) \in \bHH \times \bHH$ as $n \to \fty$. 
Our aim is to investigate conditions in which 
the sequence 
$$
\{Q(x_n,y_n)\}_{n=1}^\fty
$$ 
converges/diverges in $\D$. 
Here we say that a sequence $Q(x_n,y_n)$ {\it diverges} in $\D$ 
if it eventually exits any compact subset of $\D$, or,  
in other words, it contains no convergent subsequence. 
It is immediately follows from 
Theorem \ref{Min} (which can be viewed as a 
refinement of Thurston's Double Limit Theorem for punctured torus groups)  
that if  $(x_\fty,y_\fty)$ does not lie in $\Delta$, then 
the sequence $Q(x_n,y_n)$ converges in $\D$ 
and its limit is $Q(x_\fty,y_\fty)$. 
Thus we are interested in the case where  
$(x_\fty,y_\fty)$ lies in $\Delta$; i.e., $x_\fty=y_\fty \in \bd \HH$. 
In other words, 
we are interested in {\it exotically} convergent sequences: 

\begin{defn}\label{standexo} 
Suppose that  a sequence 
$(x_n, y_n) \in \HHD$ converges to $(x_\fty,y_\fty)$ in  $\bHH \times \bHH$,  
and that the sequence $Q(x_n,y_n)$ converges in $\D$.  
Then the convergent sequence $Q(x_n,y_n)$ is said to be {\it standard} 
if  $(x_\fty,y_\fty) \not\in \Delta$, 
and {\it exotic} if $(x_\fty,y_\fty) \in \Delta$. 
\end{defn}

We can also eliminate the case of $x_\fty=y_\fty \in \hR \sm \hQ$ 
by applying the result of Ohshika \cite{Oh1} to the case of once punctured torus:  

\begin{thm}[Ohshika]\label{Oh}
 Let $x_\fty \in \hR \sm \hQ$ and 
suppose that a sequence $(x_n,y_n) \in \HHD$ 
converges to $(x_\fty,x_\fty) \in \Delta$ as $n \to \fty$. 
Then the sequence  $Q(x_n,y_n)$ diverges in $\D$. 
\end{thm}

Thus, in what follows, we concentrate our attention to the case of  
$x_\fty=y_\fty \in \hQ$. 
By changing bases of $\pi_1(S)$ if necessary, 
we may always assume that $x_\fty=y_\fty=\fty \in \hQ$. 
Recall that $[\alpha]$ is the simple closed curve on $S$ 
corresponding to $\fty \in \hQ$, and that 
$\tau \in \Mod(S)$ denotes the Dehn twist around $[\alpha]$. 
It was first shown by 
McMullen \cite{Mc1} by using the method of 
Anderson and Canary \cite{AC} 
that there exist exotically convergent sequences: 

\begin{thm}[McMullen, Anderson-Canary]\label{AnCa}
Given $x,y \in \HH$ and an integer $p$ with $p \ne 0,-1$, 
the sequence 
\begin{equation}
\{Q(\tau^{p n} x, \tau^{(p+1)n} y)\}_{n=1}^\fty    \tag{4.1}
\end{equation}
in $\D$ converges as $n \to \fty$,
whereas the sequence 
$$
\{(\tau^{p n} x, \tau^{(p+1)n} y)\}_{n=1}^\fty
$$
 in $\HHD$ 
converges to $(\fty,\fty) \in \Delta$. 
\end{thm} 

We will see in Theorem \ref{tan} 
that every exotically convergent sequence in $\D$ is essentially 
in the form of (4.1). 
To proceed our argument, we now prepare the notion of 
horocyclical/tangential convergence for sequences in $\bHH$.   

\begin{defn}[horocyclical/tangential convergence]\label{horotan}
Let $\{x_n\} \subset \bHH$ be a convergent sequence 
with limit $x_\fty \in \bd \HH$. 
Then we say: 
\begin{itemize}
  \item   $x_n \to x_\fty$ {\it horocyclically} if for any closed horoball $B \subset \bHH$ 
touching $\bd \HH$ at $x_\fty$,   
$x_n \in B$ for all $n$ large enough. 
  \item  
  $x_n \to x_\fty$ {\it tangentially}
  if there is a closed horoball $B \subset \bHH$ 
 touching $\bd \HH$ at $x_\fty$ such that  
$x_n \not\in B$ for all $n$ large enough. 
\end{itemize}
\end{defn}

Note that the convergence 
$x_n \to \fty \in \bd \HH$ is horocyclic if and only if $\im x_n \to \fty$,  
and tangential if and only if 
$\{\im x_n\}$ are uniformly bounded above and $|\re x_n| \to \fty$. 
Therefore a sequence $\{\lam_n\} \subset \C_+$ converges 
horocyclically (resp. tangentially) to $0$  in the  sense of 
Definition \ref{horotan0} if and only if 
$\{2\pi i/\lam_n\} \subset \HH$ converges  
horocyclically (resp. tangentially) to $\fty$  in the  sense of 
Definition \ref{horotan}. 

We also remark that 
for a given $c \in \hQ$ and a sequence $\{x_n\}$ in $\HH$, 
 $x_n \to c$ horocyclically if and only if the hyperbolic lengths 
 $l_{x_n}(c)$ tend to zero, and 
  $x_n \to c$ tangentially if and only if 
 $l_{x_n}(c)$ are uniformly bounded below by a positive constant.

\subsection{Horocyclic convergence}

In this subsection, we consider the case where either $\{x_n\}$ or $\{y_n\}$ converge 
horocyclically to $\fty$.  
The case where both $\{x_n\}$ and $\{y_n\}$ converge 
tangentially to $\fty$ will be  discussed in the next subsection.  
Our main result in this subsection is the following: 

\begin{thm}\label{horo}
Suppose that a sequence $(x_n,y_n) \in \HHD$ 
converges to $(\fty,\fty) \in \Delta$ as $n \to \fty$.  
Assume that either $\{x_n\}$ or 
$\{y_n\}$ converges horocyclically to $\fty$. 
Then the sequence 
$\{Q(x_n,y_n)\}$ diverges in $\D$. 
\end{thm}

The argument in the proof of Theorem \ref{horo} is similar 
to the argument of Minsky in \cite[\S 12.3]{Mi} 
proving the properness of the map 
$x \mapsto Q(x, \fty)$ on $\bHH \sm \{\fty\}$.  
We divide the proof of Theorem \ref{horo} into the following two cases: 
\begin{description}
  \item[Case I] Both $\{x_n\}$ and $\{y_n\}$ converge  
horocyclically to $\fty$, or 
  \item[Case II] 
One of $\{x_n\}$ and $\{y_n\}$, say
  $\{x_n\}$, converges tangentially to $\fty$, and the other one,  
$\{y_n\}$,  converges horocyclically to $\fty$. 
\end{description}

Before proving Case I of Theorem \ref{horo},  
we recall the foundations of 
the thin parts of hyperbolic 3-manifolds. 
Let $M$ be a complete hyperbolic 3-manifold. 
Given $\ep>0$, the $\ep$-thin part $M^{<\ep}$ of $M$ is the set 
of all points in $M$ where the injective radius is less that $\ep$. 
The Margulis Lemma implies that 
there is a universal constant $\ep_0>0$ such that 
each component $T$ of the thin part $M^{<\ep_0}$ of $M$ 
has a standard type: either $T$ is an open solid-torus 
neighborhood of a short geodesic (called a Margulis tube), 
or $T$ is the quotient of an open horoball $B \subset \HH^3$ 
by a rank-$1$ or rank-$2$ parabolic group fixing $B$ 
(called a rank-$1$ or rank-$2$ cusp). 

\begin{proof}[Proof of  \ref{horo}: Case I] 
For simplicity, 
we assume that $x_n \ne \fty$ and $y_n \ne \fty$ for all $n$, 
but the proof for the general case is essentially the same. 
Since $\{x_n\}$, $\{y_n\}$ converge to $\fty$  horocyclically,  
both $\im x_n, \im y_n$ tend to $\fty$. 
Let $\ep_0$ be the 3-dimensional Margulis constant and let 
$\ep < \ep_0$. 
Then there is an integer $N>0$ such that $l_{x_n}([\alpha])<\ep$ and $l_{y_n}([\alpha])<\ep$ 
for all $n \ge N$. 
We fix such $n \ge N$ and denote the manifold $Q(x_n,y_n)$ by $Q_n$.  
Let $c^*$ be the geodesic realization of 
$[\alpha]$ in $Q_n$, 
and let $c^1, c^2$ be its geodesic realizations  
on each connected components of the 
boundary $\bd {\cal C}(Q_n)$ of the convex core 
${\cal C}(Q_n)$ of $Q_n$, 
where $\bd {\cal C}(Q_n)$ is regarded as a hyperbolic pleated surface. 
It is known by Sullivan (see \cite{EM}) that 
there exists a constant $K>0$,  which does not depend on $n$, 
such that the hyperbolic lengths of $c^1,c^2$ are less than $K\ep$. 

By modifying 
$c^1,c^2$ slightly in ${\cal C}(Q_n)$ in their homotopy classes if necessary, 
we may assume that both curves $c^1,c^2$ are unions of 
finitely many geodesic arcs. 
Let $A$ be an immersed annulus in 
${\cal C}(Q_n)$ with the following properties (see Figure 2): 
\begin{itemize}
\item $\bd A=c^1 \sqcup c^2$, and 
\item  $A$ is a finite union of geodesic triangles each of 
whose edges is a geodesic arc in $\bd A$ or a 
geodesic arc joining a point in $c^1$ to a point in $c^2$. 
\end{itemize}
\begin{figure}
\begin{center}
\includegraphics[width=8cm]{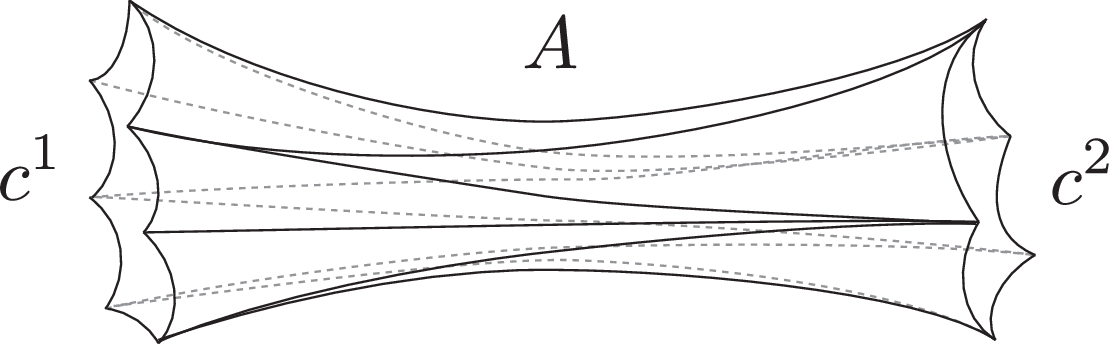}
\caption{The annulus $A$ with boundary $c^1 \sqcup c^2$. }
\label{ }
\end{center}
\end{figure}
Then one can easily check that 
$A$ is contained in the $\ep'$-thin tube ${\bf T}_{\ep'}(c^*)$ 
where $\ep'=2K\ep$. 
Let $d^*$ be the geodesic realization  in $Q_n$ of the homology class $[\beta]$  
of the generator $\beta$ of $\pi_1(S)=\la \alpha,\beta\ra$. 
Since $[\beta]$ intersects $[\alpha]$ in $S$, 
$d^*$ must intersects $A$, and hence ${\bf T}_{\ep'}(c^*)$. 

We now let $\ep \to 0$.  Then $N \to \fty$  and thus $n \to \fty$.  
Then $\ep'=2K\ep$ tends to zero and 
the distances between 
$\bd {\bf T}_{\ep'}(c^*)$ and 
$\bd {\bf T}_{\ep_0}(c^*)$ in $Q_n$ tend to $\fty$. 
Since $d^* \subset Q_n$ is not contained properly in 
${\bf T}_{\ep_0}(c^*)$ for every $n$, 
the hyperbolic lengths of 
$d^*$ in $Q_n$ diverge as $n \to \fty$.  
Since the translation lengths of the image of $\beta \in \pi_1(S)$ diverge, 
the sequence $Q(x_n,y_n)$ of representation diverges in $\D$. 
\end{proof}

Before considering Case II of Theorem \ref{horo}, 
we prepare the following  two lemmas, 
which are originally due to Kerckhoff-Thurston \cite{KT} (see also \cite{Brock}).  
 
\begin{lem}\label{triple}
Let $\{[\rho_n]\}$ be a convergent sequence in $\D$ and suppose that 
the sequence $\{\rho_n|_H:H \to \psl\}$ converges algebraically, 
where $H=\la \alpha,\,\beta^{-1}\alpha \beta \ra \subset \pi_1(S)$.
Then the sequence $\{\rho_n\}$ 
converges algebraically without taking conjugations. 
\end{lem}

\begin{proof}
Since the sequence $\{[\rho_n]\}$ converges in $\D$, 
there exist elements $\theta_n \in \psl$ such that the sequence 
$\{\theta_n \ten \rho_n \ten \theta_n^{-1}\}$ converges algebraically.   
To show that $\{\rho_n\}$ converges algebraically, 
it suffice to show that $\{\theta_n\}$ converges in $\psl$.  
Since $\{\theta_n \ten \rho_n \ten \theta_n^{-1}\}$ converges, 
$\{\theta_n \ten \rho_n|_H \ten \theta_n^{-1}\}$ converges. 
On the other hand 
$\{\rho_n|_H\}$ also converges by assumption. 
It then follows from the argument in \cite[p.35]{KT} (see also \cite{Brock}) that 
$\{\theta_n\}$ converges in $\psl$: 
In fact, let $h_1,h_2$ and $h_3$ be a triple of non-commuting elements of $H$. 
For $i=1,2$ and $3$,  $\theta_n$ maps 
the attracting fixed points $a_{n,i}$ of $\rho_n(h_i)$ to 
the attracting fixed points $a_{n,i}'$ of $(\theta_n \ten \rho_n \ten \theta_n^{-1})(h_i)$. 
Since both triples $\{a_{n,1},a_{n,2},a_{n,3}\}$ and $\{a_{n,1}',a_{n,2}',a_{n,3}'\}$ 
converge to distinct triples in $\hC$, we see that $\{\theta_n\}$ converges. 
\end{proof}

\begin{lem}\label{trick}
Let $\{[\rho_n]\}$ be a convergent sequence in $\D$ and let 
$\{\rho_n\}$ be their algebraically convergent representatives. 
Let $\{k_n\}$ be a sequence of integers. 
Then the sequence $\{\tau^{k_n} \ten [\rho_n]\}$ converges 
in $\D$  if and only if the sequence 
$\{\rho_n(\alpha)^{k_n}\}$ converges 
in $\psl$. 
\end{lem}

\begin{proof}
Let $\eta_n=\tau^{k_n} \ten \rho_n$. 
Then we have 
\begin{eqnarray*}
&&\eta_n(\alpha)=(\tau^{k_n} \ten \rho_n)(\alpha)=\rho_n(\alpha), \\
&&\eta_n(\beta)=(\tau^{k_n} \ten \rho_n)(\beta)
=\rho_n(\alpha^{k_n} \beta)
=\rho_n(\alpha)^{k_n} \rho_n(\beta). 
\end{eqnarray*}
Since $\{\rho_n\}$ converges algebraically,  
$\{\eta_n\}$ converges algebraically  if and only if 
$\{\rho_n(\alpha)^{k_n}\}$ converges in $\psl$. 

On the other hand, note  from Lemma \ref{H} that  
$\eta_n|_H=\rho_n|_H$ for all $n$.  
Since the sequence $\eta_n|_H=\rho_n|_H$ converges, 
we see from Lemma \ref{triple} 
that $\{[\eta_n]\}$ converges if and only if 
$\{\eta_n\}$ converges. 
Thus we have shown that 
$\{[\eta_n]=\tau^{k_n} \ten [\rho_n]\}$ converges if and only if 
$\{\rho_n(\alpha)^{k_n}\}$ converges. 
\end{proof}

We now back to the proof of Theorem \ref{horo}. 

\begin{proof}[Proof of  \ref{horo}: Case II]  
Recall that $\tau \in \Mod (S)$ is the Dehn twist 
around $[\alpha]=\fty \in \hQ$ 
acting on $\bHH=\ov\T(S)$ by $z \mapsto z+1$. 
To obtain a contradiction, we assume that  there exists a 
subsequence of the sequence $\{(x_n,y_n)\}_{n=1}^\fty$ (
which is also denoted by the same symbol) 
such that $\{Q(x_n,y_n)\}$ converges.  
By pass to a further subsequence if necessary, 
we may assume that there exists a divergent sequence 
$\{k_n\}$ of integers such that 
$\{\tau^{k_n}x_n\}$ converges to some 
$x_\fty' \in \bHH \sm \{\fty\}$. 
Note that $\{\tau^{k_n}y_n\}$ converges horocyclically to $\fty \in \bd \HH$ 
as well as $\{y_n\}$. 
Then we have a convergent sequence 
\begin{eqnarray*}
\tau^{k_n}  \ten Q(x_n, y_n)=Q(\tau^{k_n}x_n,\tau^{k_n}y_n) 
\to  Q(x_\fty', \fty) 
\end{eqnarray*}
in $\D$ and let 
$\rho_n \to \rho_\fty$ be a convergent sequence of their representatives.  
Since $\rho_\fty(\alpha)$ is parabolic,  the complex translation lengths 
$\lam(\rho_n(\alpha))$ of $\rho_n(\alpha)$ tend to zero as $n \to \fty$.  
Therefore, we can apply the Pivot Theorem (Theorem \ref{pivot}) to the representations 
$[\rho_n]=Q(\tau^{k_n}x_n,\tau^{k_n}y_n)$ 
to obtain 
$$
d_{\HH}\left(\frac{2\pi i}{\lam(\rho_n(\alpha))}, 
\,\tau^{k_n}x_n -\ov{\tau^{k_n}y_n}+i \right)<C
$$
for all $n$ large enough. 
Since $\tau^{k_n}x_n \to x_\fty' \ne \fty$ 
and since 
$\tau^{k_n}y_n \to \fty$ horocyclically, 
one can see that 
$$
\im \left( \frac{2\pi i}{ \lam(\rho_n(\alpha))} \right) \to \fty. 
$$
Therefore 
$\langle \rho_n(\alpha) \rangle$ converges strongly to 
$\langle \rho_\fty(\alpha) \rangle$ by Theorem \ref{cyclic} (1). 
Since we are assuming that the sequence 
$$
Q(x_n, y_n)=\tau^{-k_n} \ten [\rho_n]
$$ 
also converges in $\D$, it follows from Lemma \ref{trick} that 
$\{\rho_n(\alpha)^{-k_n}\}$ converges in $\psl$, 
which contradicts to the fact that 
$\langle \rho_n(\alpha) \rangle \to  \langle \rho_\fty(\alpha) \rangle$ strongly. 
This completes the proof. 
\end{proof}

\subsection{Tangential convergence}
In this subsection we consider the case where both $\{x_n\}$ and $\{y_n\}$ 
converge tangentially to $\fty$. 
Our main result in this case is the following:    

\begin{thm}\label{tan}
Suppose that a sequence $(x_n,y_n) \in \HHD$ 
converges to $(\fty,\fty) \in \Delta$ as $n \to \fty$, 
and that 
both  $\{x_n\}$, $\{y_n\}$ converge tangentially to $\fty$. 
Furthermore, we assume that  
there exist sequences $\{k_n\}$, $\{l_n\}$ of integers 
such that both 
$\{\tau^{k_n} x_n\}$, $\{\tau^{l_n}y_n\}$ converge 
in $\bHH \sm \{\fty\}$. 
Then the sequence $\{Q(x_n, y_n)\}$ 
converges in $\D$ if and only if 
there exist an integer $p$ such that  
\begin{eqnarray*}
(p+1)k_n-pl_n \equiv \text{const.}  
\end{eqnarray*}
for all $n$ large enough. 
\end{thm}

\begin{rem}
By passing to a subsequence if necessary, 
we may always assume that 
there exist sequences $\{k_n\}$, $\{l_n\}$ 
which satisfy the assumption in Theorem \ref{tan}.  
We also remark that the equations  $(p+1)k_n-pl_n+q \equiv 0$  $(n \gg 0)$ 
 imply  $p \ne 0,-1$ because both $|k_n|$ and $|l_n|$ tend to $\fty$. 
\end{rem}

One of the key observations in the proof of Theorem \ref{tan} is the following: 

\begin{thm}\label{prime}
Let $u_n \to u_\fty$, $v_n \to v_\fty$ be 
convergent sequences in $\bHH \sm \{\fty\}$, 
and let $\{m_n\}$ be a sequence of integers such that $|m_n| \to \fty$.  
Then 
$Q(u_n, \tau^{m_n}v_n) \to Q(u_\fty,\fty)$ in $\D$, and 
let $\rho_n \to \rho_\fty$ be a convergent sequence of their representatives. 
In this situation, we have the following: 
\begin{enumerate}
  \item  
  $\lim_{n \to \fty}\rho_n(\alpha)^{m_n}$ exists. 
\item
 The sequence 
$\la \rho_n(\alpha) \ra$  of cyclic groups  converges 
geometrically to a rank-$2$ parabolic group 
$\la \delta,\hat \delta \ra$ 
where $\delta=\rho_\fty(\alpha)$ and 
 $\hat \delta=\lim_{n \to \fty}\rho_n(\alpha)^{m_n}$. 
\end{enumerate}
\end{thm}

\begin{rem}
The essential points in Theorem \ref{prime} are  that 
the sequence $\la \rho_n(\alpha) \ra$ converges 
geometrically without passing to a subsequence, 
and that the limit 
 $\hat \delta=\lim_{n \to \fty}\rho_n(\alpha)^{m_n}$ 
is primitive in the geometric limit. 
\end{rem}

\begin{proof}
(1) 
By changing the marking of the convergent sequence
\begin{eqnarray*}
Q(u_n, \tau^{m_n}v_n) \to Q(u_\fty,\fty), 
\end{eqnarray*}
we obtain another  convergent sequence
$$
\tau^{-m_n} \ten Q(u_n, \tau^{m_n}v_n) 
= Q(\tau^{-m_n}u_n, v_n) \to Q(\fty,v_\fty). 
$$ 
We let 
$$
\bar{\rho}_n=\tau^{-m_n} \ten \rho_n 
$$
so that $\bar{\rho}_n$ represent $Q(\tau^{-m_n}u_n, v_n)$. 

We first  show that $\{\bar \rho_n\}$ converges algebraically.  
Since $\{\rho_n\}$ converges,  and since 
 $\bar \rho_n|_H=\rho_n|_H$ for all $n$ by Lemma \ref{H}, 
$\{\bar \rho_n|_H\}$ also converges. 
In addition, since $\{[\bar{\rho}_n]\}$ converges in $\D$, 
we see from Lemma \ref{triple} that 
$\{\bar \rho_n\}$ converges algebraically.  
Note that  the limit $\bar{\rho}_\fty=\lim_{n \to \fty}\bar \rho_n$ 
of the sequence $\{\bar \rho_n\}$ is 
a representative of $Q(\fty,v_\fty)$. 

Now one see from 
$\bar\rho_n=\tau^{-m_n} \ten \rho_n$ that 
$$
\bar\rho_n(\beta)=(\tau^{-m_n} \ten \rho_n) (\beta)
=\rho_n(\alpha^{-m_n} \beta)=\rho_n(\alpha)^{-m_n} \rho_n(\beta), 
$$
which is rewritten as 
$\rho_n(\alpha)^{m_n}=\rho_n(\beta) \bar\rho_n(\beta)^{-1}$. 
Therefore $\{\rho_n(\alpha)^{m_n}\}$ has a limit: 
$$
\lim_{n \to \fty} 
  \rho_n(\alpha)^{m_n}
=\rho_\fty(\beta) \bar \rho_\fty(\beta)^{-1}. 
$$

(2) 
To show that $\la \rho_n (\alpha)\ra \to \la \delta,\hat \delta\ra$ geometrically, 
it suffices to show that any subsequence $\{n'\}$ of $\{n\}_{n=1}^\fty$ has a further 
subsequence $\{n''\}$ such that 
$\la \rho_{n''}(\alpha)\ra \to \la \delta,\hat \delta\ra$ geometrically.  
Applying the Pivot Theorem (Theorem \ref{pivot}) to the representations 
$[\rho_{n'}]=Q(u_{n'},\tau^{m_{n'}}v_{n'})$, we obtain 
$$
d_{\HH} \left( \frac{2\pi i}{\lam(\rho_{n'}(\alpha))}, u_{n'}-\ov{\tau^{m_{n'}}v_{n'}} +i  \right)
=
d_{\HH} \left( \frac{2\pi i}{\lam(\rho_{n'}(\alpha))}, (u_{n'}-\ov{v_{n'}}+i)-m_{n'}  \right)
<C
$$
for all ${n'}$ large enough. 
Since $\{u_{n'}-\ov{v_{n'}} +i\}$ has a limit $u_\fty-\ov{v_\fty}+i$ in $\HH$,  
we can take a subsequence $\{n''\}$of $\{n'\}$ so that the sequence
$$
\frac{2 \pi i}{\lam(\rho_{n''}(\alpha))} +m_{n''}
$$ 
has a limit in $\HH$. 
It then follows from Theorem \ref{cyclic} (2) that 
$\la \rho_{n''}(\alpha) \ra \to \la \delta'', \hat \delta'' \ra$ geometrically, 
where $\delta''=\lim_{n'' \to \fty} \rho_{n''}(\alpha)$ and 
$\hat \delta''=\lim_{{n''} \to \fty} \rho_{n''}(\alpha)^{m_{n''}}$. 
Since $\delta''=\delta$ and $\hat \delta''=\hat \delta$, we obtain 
$\la \rho_{n''}(\alpha) \ra \to \la \delta, \hat \delta \ra$ geometrically. 
Therefore $\la \rho_{n}(\alpha) \ra \to \la \delta, \hat \delta \ra$ geometrically 
without passing to  a subsequence. 
\end{proof}

We now back to the proof of Theorem \ref{tan}. 

\begin{proof}[Proof of Theorem \ref{tan}]
We begin by fixing our notation. 
Let 
$$u_n:=\tau^{k_n}x_n, \quad v_n:=\tau^{l_n}y_n 
$$ 
and let $u_\fty, v_\fty \in \bHH \sm \{\fty\}$ denote  
the limits of $\{u_n\}$, $\{v_n\}$, respectively. 
 Let 
$$
m_n:=k_n-l_n
$$ 
and let us consider the re-marked sequence 
$$
\tau^{k_n} \ten Q(x_n, y_n)=Q(\tau^{k_n} x_n, \tau^{k_n} y_n)=Q(u_n, \tau^{m_n} v_n)
$$
of $Q(x_n,y_n)$. 

To show Theorem \ref{tan}, 
we first show that the following two conditions are equivalent 
under the assumption that $|m_n|=|k_n-l_n| \to \fty$: 
\begin{enumerate}
\item $\{Q(x_n,y_n)\}$ converges. 
\item  There exist an integer $p$ such that $(p+1)k_n-pl_n \equiv \text{const.}$
 for all $n$ large enough. 
\end{enumerate}
After that, we will show that each of both conditions (1) and (2) 
implies $|m_n|\to \fty$ to complete the proof.  

We now assume that  $|m_n| \to \fty$. Then 
$$
Q(u_n,\tau^{m_n} v_n) \to Q(u_\fty, \fty)
$$
in $\D$,  and 
let $\rho_n \to \rho_\fty$ be a convergent sequence of their representatives. 
It then follows from  Theorem \ref{prime} that   
$$
\la \rho_n(\alpha) \ra \to \la \delta, \hat\delta\ra
$$ geometrically,  
where $\delta=\rho_\fty(\alpha)$ and 
$\hat\delta=\lim_{n \to \fty}\rho_n(\alpha)^{m_n}$. 
Recall from Lemma \ref{trick} that the sequence 
$$
Q(x_n, y_n) =\tau^{-k_n} \ten Q(u_n, \tau^{m_n} v_n)=\tau^{-k_n} \ten [\rho_n]
$$
converges if and only if 
$\{\rho_n(\alpha)^{-k_n}\}$ converges in $\psl$.   
Therefore we need to show that $\{\rho_n(\alpha)^{-k_n}\}$ converges 
if and only if the condition (2) is satisfied. 

Now suppose that  $\{\rho_n(\alpha)^{-k_n}\}$ converges. 
Then its limit must lie in $\la \delta, \hat\delta\ra$, and thus 
is of the form 
${\hat \delta}^p \delta^q$  
for some integers $p,\,q$.  
Since $\{\rho_n(\alpha)^{pm_n+q}\}$ also converges to ${\hat \delta}^p \delta^q$, 
it follows from  discreetness 
that  $\rho_n(\alpha)^{-k_n} \equiv \rho_n(\alpha)^{pm_n+q}$ 
for all $n$ large enough (see Lemma 3.6 in \cite{JM}). 
Therefore  $-k_n \equiv p  m_n+q \,(n \gg 0)$.   
Using $m_n=k_n-l_n$,  we obtain the equations 
$(p+1) k_n-p l_n \equiv \text{const.}$ for all $n$ large enough. 
Conversely, suppose that there exist integers $p,\,q$ such that 
$(p+1) k_n-p l_n+q \equiv 0$ for all $n$ large enough. 
Using $m_n=k_n-l_n$ we have 
$\rho_n(\alpha)^{-k_n}=\rho_n(\alpha)^{pm_n+q} \to \hat \delta^p \delta^q$. 
Thus $\{\rho_n(\alpha)^{-k_n}\}$ converges. 
So far we have shown the theorem under the assumption that $|m_n| \to \fty$. 

We now show that the convergence of $\{Q(x_n,y_n)\}$ implies $|m_n| \to \fty$. 
Suppose that the sequence 
$Q(x_n,y_n)$ converges in $\D$. 
In addition, suppose for contradiction that there exists a subsequence 
of $\{m_n\}$ (which is also denoted by $\{m_n\}$) 
and an integer $m'$ 
such that $m_n \equiv  m'$  for all $n$ large enough.   
Then  $Q(u_n, \tau^{m_n} v_n) \to Q(u_\fty, \tau^{m'} v_\fty)$ in $\D$. 
We let $\rho_n \to \rho_\fty$ be a convergent sequence of their representatives. 
Since both $u_\fty$ and $\tau^{m'} v_\fty$ lie in $\bHH \sm \{\fty\}$, 
$\rho_\fty(\alpha)$ is loxodromic. 
Therefore $\la \rho_n(\alpha) \ra \to \la \rho_\fty(\alpha) \ra$  strongly. 
Since $|k_n| \to \fty$,  the sequence 
$\rho_n(\alpha)^{k_n}$ diverges in $\psl$. 
It then follows from Lemma \ref{trick}  that the sequence 
$Q(x_n, y_n) =\tau^{k_n} \ten Q(u_n, \tau^{m_n} v_n)$ diverges in $\D$, 
which is a contradiction. Thus we obtain $|m_n| \to \fty$. 

Finally, we show that the condition 
$(p+1) k_n - p l_n \equiv \text{const.} \, (n \gg 0)$ implies $|m_n| \to \fty$. 
This follows immediately from $(p+1)k_n-pl_n=k_n+pm_n$ and $|k_n| \to \fty$. 
\end{proof}

In the next section, we will give an explicit description of the limits of exotically 
convergent sequences.

\section{Bers and Maskit slices and the Maskit embedding} 

In this section, we recall Bers and Maskit slices of $\D$, 
and a natural embedding of a Maskit slice into 
the complex plane $\C$ which is called the Maskit embedding. 
We will rephrase  Theorem \ref{tan} in terms of the Maskit embedding.
We will also explain Anderson-Canary's wrapping construction. 

\subsection{Bers and Maskit slices}

Given $y \in \bHH$, we define an embedding 
\begin{eqnarray*}
b_y:
\begin{cases}
\bHH \to \D    & \text{if } y \in \HH \\
\bHH \sm \{y\}  \to \D   & \text{if } y \in \bd \HH 
\end{cases}
\end{eqnarray*}
by $x \mapsto Q(x,y)$,  
and denote the image of $b_y$ by 
$$
\B_y:=\{Q(x,y)\,|\,x \in \bHH, \, (x,y) \not\in \Delta\}.  
$$
The image $\B_y$ is called a {\it Bers slice} if $y \in \HH$, 
and a {\it Maskit slice} if $y \in \hQ \subset \bd \HH$.  
However, if there is no confusion, 
we generally called the set $\B_y$ for arbitrary $y \in \bHH$ as a {\it Bers slice} 
for simplicity. 
By Theorem \ref{Min}  the map 
$b_y:\bHH \to \B_y$ 
is a homeomorphism  if  $y \in \HH$. 
It is also known that  the map 
$b_y:\bHH \sm \{y\}  \to \B_y$ is a homeomorphism 
even when $y \in \bd \HH$
 by Theorem \ref{Min} together with the arguments 
by Minsky \cite{Mi} in the case of $y \in \hQ$,  
and by Ohshika \cite{Oh1} in the case of $y \in \hR \sm \hQ$. 

Given $x \in \bHH$, we similarly define a slice $\B_x^* \subset \D$ by 
$$
\B_x^*=\{Q(x,y)\,|\,y \in \bHH,  \,  (x,y) \not\in \Delta \}. 
$$

\subsection{The Maskit embedding}

We now explain a natural embedding of 
the Maskit slice 
$\B_\fty$ for $\fty \in \hQ$ into the complex plane $\C$, 
which is known as the Maskit embedding. 
We refer the reader to \cite{KS}  for more details. 

Given $z, \mu \in \C$, we define elements $T_z, \,U_\mu$ of $\psl$ by 
\begin{eqnarray*}
T_z 
=\left(
  \begin{array}{cc}
    1   &  z  \\
    0   &  1  \\
  \end{array}
\right), \quad 
U_\mu=\left(
  \begin{array}{cc}
    i\mu   &  i  \\
    i   & 0   \\
  \end{array}
\right). 
\end{eqnarray*}
We will later use the fact that 
$U_\mu U_\nu^{-1}=T_{\mu-\nu}$ holds for any $\mu,\nu \in \C$. 
Given $[\rho] \in \RR(S)$ such that $\rho(\alpha)$ is parabolic,  
there is a unique $\mu \in \C$ such that 
$\rho$ is conjugate to a representation 
$\brho_\mu:\pi_1(S) \to \psl$ defined by 
\begin{eqnarray*}
\brho_\mu(\alpha)=T_2=\left(
  \begin{array}{cc}
    1   &  2 \\
    0   &  1  \\
  \end{array}
\right), \quad 
\brho_\mu(\beta)=U_\mu=\left(
  \begin{array}{cc}
    i\mu   &  i  \\
    i   & 0   \\
  \end{array}
\right).  
\end{eqnarray*}
The representation $\brho_\mu$ is  normalized so that the fixed points of 
$\brho_\mu(\alpha)$, $\brho_\mu(\beta^{-1} \alpha \beta)$ and 
$\brho_\mu(\alpha^{-1} \beta^{-1} \alpha \beta)$ coincide with $\fty$, $0$ and $-1$, respectively. 
One can see that the map 
$$
\Phi_\fty:\C \to \RR(S), \quad  \mu \mapsto [\brho_\mu] 
$$
is a holomorphic embedding,  and 
from Theorem \ref{Min} that 
\begin{eqnarray*}
\Phi_\fty(\C) \cap \D&=&\B_\fty \sqcup \B_\fty^*.
\end{eqnarray*}
We let denote the preimages of $\B_\fty$ and $\B_\fty^*$ in $\C$ by 
\begin{eqnarray*}
\M=\Phi_\fty^{-1}(\B_\fty) \quad \text{and}\quad \M^*= \Phi_\fty^{-1}(\B_\fty^*). 
\end{eqnarray*}
It is known that the subset $\M \subset \C$ is contained in the upper half plane $\HH$,  
and that $\M^*=\{\ov{\mu}\,|\,\mu \in \M\}$ is the complex conjugation of $\M$. 
The set $\M$ is also called the Maskit slice. 
Now we define  a bijective map 
$$
\Psi:\bHH \sm \{\fty\} \to \M 
$$
by $\Psi=\Phi_\fty^{-1} \circ b_\fty$. 
Then we have 
$[\brho_{\Psi(x)}]=Q(x,\fty)$ and 
$[\brho_{\ov{\Psi(x)}}]=Q(\fty,x)$ for every $x \in \bHH \sm \{\fty\}$. 

\subsection{Limits of exotically  convergent sequences}

The next theorem gives us a precise description of the limit of an 
exotically convergent sequence in Theorem \ref{tan}. 

\begin{thm}\label{limit}
Suppose that a sequence $(x_n,y_n) \in \HHD$ 
satisfies the same assumption as in Theorem \ref{tan}. 
Let 
$u_\fty, \,v_\fty \in \bHH \sm \{\fty\}$ be  
the limits of the sequences 
$u_n=\tau^{k_n} x_n$, $v_n=\tau^{l_n}y_n$, and let 
$\mu,\, \nu \in \M$ be the corresponding points of $u_\fty,\,v_\fty$ 
via the map $\Psi:\bHH \sm \{\fty\} \to \M$, respectively. 
Assume that 
there exist integers $p,\,q$ 
such that $(p+1)k_n-p l_n+q \equiv 0$ for all $n$ large enough. 
Then we have 
\begin{eqnarray*}
 \lim_{n \to \fty}Q(x_n,y_n)  = [\brho_\xi] 
\end{eqnarray*}
with
$$
\xi=(p+1)\mu-p \ov{\nu}+2q. 
$$
\end{thm}

\begin{proof}
Let $m_n:=k_n-l_n$. 
Then the equation $(p+1)k_n-p l_n+q=0$ 
can be rewritten as $k_n=-p m_n-q$. 
As observed in the proof of Theorem \ref{tan} 
we have $|m_n|=|k_n-l_n| \to \fty$. 
Thus the sequence $Q(u_n,\tau^{m_n}v_n)$ 
converges to $Q(u_\fty,\fty)=[\brho_\mu]$ as $n \to \fty$.  
Choose representatives $\rho_n$ of  $Q(u_n,\tau^{m_n}v_n)$ 
so that $\rho_n \to \brho_\mu$ algebraically. 
Then $\delta:=\lim_{n \to \fty}\rho_n(\alpha)=\brho_\mu(\alpha)=T_2$. 

We next show that 
$\hat \delta:=\lim_{n \to \fty}\rho_n(\alpha)^{m_n}=T_{\mu-\ov\nu}$. 
Let $\bar{\rho}_n:=\tau^{-m_n} \ten \rho_n$ 
so that $\bar{\rho}_n$ represent $Q(\tau^{-m_n}u_n, v_n)$. 
Then as observed in the proof of Theorem \ref{tan}, 
$\{\bar\rho_n\}$ converges algebraically, and the limit 
$\lim_{n \to \fty} \bar\rho_n$ represent 
$Q(\fty,v_\fty)=[\brho_{\ov \nu}]$. 
Since $\rho_n|_H=\bar\rho_n|_H$, $\brho_\mu|_H=\brho_{\ov \nu}|_H$
and since $\rho_n \to \brho_\mu$, 
we see that $\bar\rho_n|_H \to \brho_{\ov \nu}|_H$. 
In addition, since $[\bar\rho_n] \to [\brho_{\ov \nu}]$, 
it follows from Lemma \ref{H} that $\bar\rho_n \to \brho_{\ov \nu}$. 
Therefore one see form the proof of Theorem \ref{prime} that 
$$
\hat \delta:=\lim_{n \to \fty}\rho_n(\alpha)^{m_n}=
\brho_{\mu}(\beta)\brho_{\ov{\nu}}(\beta)^{-1}=
U_\mu U_{\ov{\nu}}^{-1}=T_{\mu-\ov\nu}. 
$$

Now let $\eta_n:=\tau^{-k_n}\ten\rho_n$ so that it represent 
$\tau^{-k_n} \ten Q(u_n,\tau^{m_n} v_n)=Q(x_n,y_n)$. 
One can observe that $\eta_n \to \brho_\xi$ with $\xi=(p+1)\mu-p \ov{\nu}+2q$ 
as follows: 
\begin{eqnarray*}
\eta_n(\alpha)&=&\tau^{-k_n} \ten \rho_n(\alpha)=\rho_n(\alpha) \\
& \to& \brho_\mu(\alpha)=T_2, \\ 
\eta_n(\beta)
&=&\tau^{-k_n} \ten \rho_n(\beta)
= \rho_n(\alpha^{-k_n}\beta)
=\rho_n(\alpha)^{pm_n+q}\rho_n(\beta) \\
&\to& \hat\delta^p \delta^q \brho_\mu(\beta)
=(T_{\mu-\ov\nu})^p(T_2)^q U_\mu =U_\xi. 
\end{eqnarray*} 
Therefore we have shown that the sequence 
$[\eta_n]=Q(x_n,y_n)$ converges to $[\brho_\xi]$. 
\end{proof}

\subsection{ Anderson-Canary's wrapping construction}
Here we remark connection between our result and Anderson-Canary's 
wrapping construction. 
More information can be found in  \cite{AC}, \cite{Brom} and \cite{Mc1}. 
 
Let $\mu, \nu \in \M$.  
Then the group 
$$
\hat \Gamma=\la \brho_\mu(\pi_1(S)), \brho_{\ov{\nu}}(\pi_1(S))\ra 
=\la \brho_\mu(\alpha), \brho_\mu(\beta), \brho_{\ov \nu}(\beta)\ra
=\la T_2,U_\mu,U_{\ov \nu}\ra
$$
generated by the images of $\brho_\mu$ and $\brho_{\ov \nu}$ uniformizes
a manifold $M_{\hat \Gamma}=\HH^3/\hat \Gamma$ 
homeomorphic to $S \times (0,1)$ with simple closed curve $\alpha \times \{1/2\}$ removed: 
$$
M_{\hat \Gamma} \cong S \times (0,1) \sm \alpha \times \{1/2\}. 
$$
Here the subgroup of $\pi_1(M_{\hat \Gamma}) \cong \hat \Gamma$ carried by the surface 
$S \times \{1/4\} \subset M_{\hat \Gamma}$ is equal to 
$\brho_\mu(\pi_1(S))$ and  
the subgroup  carried by the surface $S \times \{3/4\}$ is equal to 
$\brho_{\ov\nu}(\pi_1(S))$.  
The rank-2 parabolic subgroup 
$$
\la \brho_\mu(\alpha), \brho_\mu(\beta) \brho_{\ov \nu}(\beta)^{-1}\ra=\la T_2, T_{\mu-\ov{\nu}}\ra
$$
of $\hat \Gamma$ corresponds to the fundamental group of the rank-2 cusp of 
$M_{\hat \Gamma}$.  

For a given integer $p$,  we define a {\it wrapping map}  
$$
w_p:S \to M_{\hat \Gamma}
$$
as follows (see Figure 3): 
\begin{enumerate}
\item The wrapping map $w_p$  is an immersion determined up to homotopy. 
\item When $p=0$, we let 
$w_0:S \to S \times \{1/4\} \subset M_{\hat \Gamma}$ be  the identity map.  
\item For general $p$, we let $w_p:S \to M_{\hat \Gamma}$ be an immersion such that the image 
$w_p(S)$ in $M_{\hat \Gamma}$
is obtained by cutting $S \times \{1/4\}$ along $\alpha \times \{1/4\}$ and inserting 
 at the cut locus 
an annulus which wraps $p$ times around the rank-2 cusp 
$\alpha \times \{1/2\}$ of $M_{\hat \Gamma}$.
It is also required that $w_p$ 
is homotopic to $w_0$ in $S \times (0,1)$. 
\end{enumerate}
\begin{figure}
\begin{center}
\includegraphics[width=8cm]{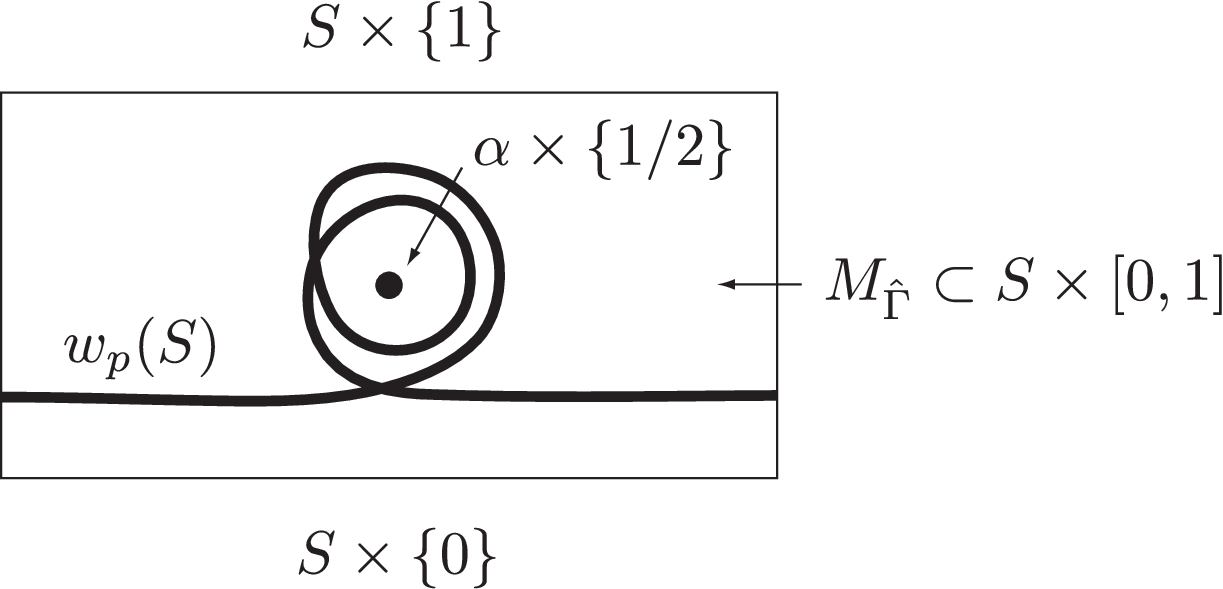}
\caption{ Schematic figure of the image $w_p(S)$ of the wrapping map $w_p$ ($p=2$ case). }
\label{ }
\end{center}
\end{figure}
We remark that when $p=-1$, 
$w_{-1}$ is homotopic to the identity map 
$S \to S \times \{3/4\}$ in $M_{\hat \Gamma}$. 

The group homomorphism 
$(w_0)_*:\pi_1(S) \to \psl$ induced by $w_0$ is conjugate to $\brho_\mu$. 
Similarly $(w_{-1})_*$ is conjugate to $\brho_{\ov{\nu}}$. 
In general,  one see that $(w_p)_*$ is conjugate to the representation 
\begin{eqnarray*}
&&\alpha \mapsto \brho_\mu(\alpha)=T_2, \\
&&\beta \mapsto \brho_\mu(\beta)(\brho_\mu(\beta) \brho_{\ov \nu}(\beta)^{-1})^p=T_\mu (T_{\mu-\ov\nu})^p. 
\end{eqnarray*}
Therefore $(w_p)_*$ is conjugate to 
$\brho_{\xi}$ with  
$$
\xi=(p+1)\mu-p \ov\nu. 
$$  
Since $\brho_\xi(\pi_1(S)) \subset \hat \Gamma$ is discrete, 
one see that $\xi \in \M$ if $p \ge 0$ and $\xi \in \M^*$ if $p \le -1$. 

Now let 
$$
F_n: M_{\hat \Gamma} \to M_n
$$ 
be the  $(1,n)$ Dehn filling of  $M_{\hat \Gamma}$ at the rank-2 cusp  
without changing the end invariants $x,y \in \ov{\T}(S)=\ov\HH$ of $M_{\hat \Gamma}$. 
Then we obtain a sequence 
$$
[(F_n \circ w_p)_*]=Q(\tau^{pn}x,\tau^{(p+1)n}y)
$$ 
in $\D$  converging to $[\brho_\xi]$ as $n \to \fty$, see \cite{Mc1}. 
This is Anderson-Canary's wrapping construction of exotically convergent sequence. 
Note that in this case $k_n:=-pn$ and $l_n:=-(p+1)n$ satisfy $(p+1)k_n-pl_n \equiv 0$. 
One can easily see from Theorems \ref{tan} and \ref{limit} that every exotically convergent sequence 
is obtained by Anderson-Canary's wrapping construction. 

\subsection{Subsets of the Maskit slice} 

For later use, here we study some properties of the set of $\xi \in \C$ in Theorem \ref{limit}. 

For an given integer $p$, we define 
a subset $\M(p)$ of $\C$ by 
\begin{eqnarray*}
\M(p)&=&\{(p+1) \mu-p \ov{\nu} \in \C\,|\, \mu,\nu \in \M\}. 
\end{eqnarray*}
Note that $\M(0)=\M$,  and that for every $p$,  
the set $\M(p)$ is invariant under the translation 
$z \mapsto z+2$. 
The set $\M(p)$ is equal to the set of points $\xi \in \C$ such that $[\brho_\xi]$ are induced by
wrapping maps $w_p$. 

We denote by $\M^*(p)$ the complex conjugation of $\M(p)$: 
$$
\M^*(p)=\{\ov\xi \,|\,\xi \in \M(p)\}. 
$$ 
Then one can easily check that 
$$
\M^*(p)=\M(-p-1) 
$$ 
for every $p$. 
For this reason, we mainly use the notation $\M(p), \M^*(p)$ for non-negative integers $p \ge 0$. 
We prepare some basic properties of $\M(p)$ in the following lemma. 

\begin{lem}\label{subset}
\begin{enumerate}
\item $\M(p) \subset \M$ for every integer $p \ge 1$. 
\item $\M(p) \subset \M(1)$ for every integer $p \ge 2$. 
\end{enumerate}
\end{lem}

\begin{proof}
(1) 
For any $\mu,\nu \in \M$ and an integer $p \ge 1$, 
we see from the argument in section 5.4 that 
$(p+1) \mu-p \ov{\nu} \in \M$. 
This implies that $\M(p) \subset \M$. 

(2)
Suppose that $p \ge 2$. 
For any $\mu,\nu \in \M$, we want to show that 
$(p+1)\mu-p\ov{\nu}$ lies in $\M(1)$. 
Let  $k \ge 1$ so that $p=k+1$. 
Then we have 
$$
(p+1)\mu-p\ov{\nu}=2\mu-\ov{(k+1)\nu-k\ov{\mu}}. 
$$
Let $\xi:=(k+1)\nu-k\ov{\mu}$. Then 
$\xi$ lies in $\M(k)$, and hence in $\M$ from (1).   
Therefore $(p+1)\mu-p\ov{\nu}=2\mu-\ov{\xi}$ lies in $\M(1)$. 
\end{proof}

\section{Hausdorff limits of Bers slices}

In this section, we consider Hausdorff limits of sequences of Bers slices. 
Especially, we obtain a condition such that the Hausdorff limit is strictly 
bigger than a Bers slice. 

Let $y_n \to y_\fty$ be a convergent sequence in $\bHH$. 
Then by Theorem \ref{Min} the maps 
\begin{eqnarray*}
b_{y_n}:\{x \in \bHH  \,|\,(x,y_n) \not\in \Delta\} \to \B_{y_n},  
\quad x \mapsto Q(x,y_n)
\end{eqnarray*}
converge to the map 
\begin{eqnarray*}
b_{y_\fty}:\{x \in \bHH  \,|\, (x,y_\fty) \not\in \Delta\} \to \B_{y_\fty}, 
\quad x \mapsto Q(x,y_\fty)
\end{eqnarray*} 
at each point in $\{x \in \bHH  \,|\, (x,y_\fty) \not\in \Delta\}$. 
After passing to a subsequence if necessary, 
the sequence  $\{\B_{y_n}\}_{n=1}^\fty$ 
converges to a closed subset  
$\hat\B \subset \D$ in the Hausdorff topology on $\D$; that is, 
\begin{enumerate}
\item for any $[\rho] \in  \hat\B$ there exists a sequence $[\rho_n] \in \B_{y_n}$ 
such that $[\rho_n] \to [\rho]$, and 
\item if elements $[\rho_{n_j}] \in \B_{y_{n_j}}$ converge to 
$[\rho]$,  then $[\rho] \in \hat\B$. 
\end{enumerate}
One can easily check that the Bers slice $\B_{y_\fty}$ 
for $y_\fty=\lim_{n \to \fty} y_n$ is contained in the Hausdorff limit  $\hat\B$.  
In the next theorem, we consider the situation that 
$\B_{y_\fty}$ is properly contained in the Hausdorff limit  $\hat\B$.   

\begin{thm}\label{bersgeom}
Let $y_n \to y_\fty$ be a convergent sequence in $\bHH$, 
and suppose that 
$\{\B_{y_n}\}_{n=1}^\fty$ converges to $\hat\B$ in the sense of Hausdorff.  
Then we have the following: 
\begin{enumerate}
\item 
$\B_{y_\fty} \subsetneq \hat\B$ 
if and only if $y_\fty \in \hQ$ and 
$y_n \to y_\fty$ tangentially.  
\item Suppose that $y_n \to y_\fty$ tangentially and that $y_\fty=\fty \in \hQ$. 
We further assume that there exists a sequence $\{l_n\}$ 
of integers such that $\{\tau^{l_n} y_n\}$ converges to some 
$v_\fty \in \bHH \sm \{ \fty\}$. 
Let $\nu \in \M$ be the corresponding point of $v_\fty$ via the map 
$\Psi:\bHH \sm\{\fty\} \to \M$. 
Then we have 
\begin{eqnarray*}
 \hat\B=\{[\brho_\xi]  \,|\, \xi \in \M \sqcup (\M^*+2\ov{\nu})\}. 
\end{eqnarray*}
In particular 
$\B_{\fty} \subsetneq \hat\B \subsetneq \B_{\fty} \sqcup \B_{\fty}^*$.  
\end{enumerate}
\end{thm}

\begin{proof}
(1) 
Assume first that 
$\B_{y_\fty}\subsetneq \hat\B$. 
Then there exists a sequence $\{x_n\}$ in $\bHH$ such that the sequence 
$Q(x_n,y_n) \in \B_{y_n}$ converges in $\D$ to a point which does not lie in 
$\B_{y_\fty}$. 
By pass to a subsequence if necessary, we may assume  that  
$\{x_n\}$ converges to some $x_\fty \in \bHH$.  
If $(x_\fty,y_\fty) \not\in \Delta$, 
Theorem \ref{Min} implies that 
$Q(x_n,y_n) \to Q(x_\fty, y_\fty) \in \B_{y_\fty}$, 
which contradicts our assumption. 
Therefore we have $(x_\fty,y_\fty) \in \Delta$, i.e.,  $x_\fty=y_\fty \in \bd \HH$. 
Since the sequence 
$Q(x_n, y_n)$ converges in $\D$, we see from  Theorem \ref{Oh} that 
$y_\fty$ must lie in $\hQ \subset \bd \HH$, 
and  from Theorem \ref{horo} that $\{y_n\}$ converges  to $ y_\fty$ tangentially.  

The converse is deduced from the statement of (2). 

(2)  
We first show  for a given sequence $\{Q(x_n,y_n) \in\B_{y_n}\}_{n=1}^\fty$ that 
every accumulation points of this sequence  must 
lie in the set 
\begin{equation}
\{[\brho_\xi]\,|\,\xi \in \M \sqcup (\M^*+2\bar\nu)\}.  \tag{6.1}
\end{equation}
In fact, pass to a subsequence so that the sequence 
$Q(x_n,y_n)$ converges to some $[\rho_\fty] \in \D$,  
and that $\{x_n\}$ converges to some $x_\fty \in \bHH$.  
If $x_\fty$ does not equal to 
$y_\fty=\fty$ then $[\rho_\fty]=Q(x_\fty, \fty)$ by Theorem \ref{Min}. 
In this case we have $[\rho_\fty] \in \{[\brho_\xi]\,|\,\xi \in \M\}$.  
Now let us consider the case that $x_\fty=y_\fty=\fty$. 
Since $Q(x_n,y_n)$ converges, 
one see from Theorem \ref{horo} that $\{x_n\}$ converges tangentially to $\fty$ 
as well as $\{y_n\}$.  
Recall that we are assuming that 
there is a sequence $\{l_n\}$ of integers such that the sequence $v_n=\tau^{l_n}y_n$ converges  
to $v_\fty  \in \bHH \sm \{\fty\}$. 
Passing to a further subsequence,  
we may also assume that 
there is a sequence $\{k_n\}$ of integers  
such that the sequence $u_n=\tau^{k_n} x_n$ converges to some 
$u_\fty  \in \bHH \sm \{\fty\}$. 
Let $\mu=\Psi(u_\fty),\, \nu=\Psi(v_\fty)$. 
Then by Theorem \ref{limit} 
the sequence $Q(x_n,y_n)$ converges to $[\rho_\fty]=[\brho_\xi]$ 
with $\xi=(p+1)\mu-p\ov\nu+2q$ for some integers $p,q$. 
In addition $p\ne 0, -1$ because  $Q(x_n,y_n)$ is an exotically  convergent sequence. 
If $p \ge 1$ then $\xi \in \M$. 
If $p \le -2$, we have 
$$
\xi=(p+1)\mu-p\ov\nu+2q=\ov{(k+1)\nu-k \ov{\mu}+2q}+2 \ov{\nu}
$$
where $k=-p-1 \ge 1$, and hence $\xi \in \M^*+2\ov\nu$ from Lemma \ref{subset}. 
In both cases the limit $[\rho_\fty]=[\brho_\xi]$ of the sequence $Q(x_n,y_n)$ must lies in the set (6.1). 

We next show that 
for a given  $\xi \in \M \sqcup (\M^*+2\bar\nu)$, 
$[\brho_\xi]$ is the limit of a sequence $\{Q(x_n,y_n) \in \B_{y_n}\}_{n=1}^\fty$. 
When $\xi \in \M$, let $x=\Psi^{-1}(\xi)$. 
Then the sequence $Q(x,y_n) \in \B_{y_n}$  converges to $Q(x,\fty)=[\brho_\xi]$. 
When $\xi \in \M^*+2\ov\nu$, take $\mu \in \M$ such that $\xi=-\mu+2\ov\nu$. 
Let $u=\Psi^{-1}(\mu)$ 
and $x_n=\tau^{-2l_n}u$. 
Then the sequence $(x_n,y_n) \in \HHD$ satisfies the condition 
of Theorem \ref{limit} in the case of $p=-2$ and $q=0$, 
and thus the sequence 
$$
Q(x_n,y_n)=Q(\tau^{-2l_n}u, \tau^{-l_n}v_n) \in \B_{y_n}
$$
converges to $[\brho_\xi]$ with $\xi=-\mu+2\ov\nu$.  

Thus we have shown that  the geometric limit $\hat\B$ equals the set  (6.1). 
The last sentence is obvious. 
\end{proof}

\section{Self-Bumping of $\D$}

In this section, 
we consider sequences which give rise to the self-bumping of $\D$, and  
obtain a precise description of the 
set of points at which $\D$ self-bumps. 

To state our result, we need to define subsets of Maskit slices $\B_y$ for $y \in \hQ$  
associated to the subsets $\M(p)$ of $\M$.
For every $y \in \hQ$, 
let $\sigma_y \in \Mod(S)$ denote a mapping class which takes 
the homology class  $\pm [\alpha] \in H_1(S)$ associated to $\fty \in \hQ$ 
to  the homology class $\pm(s [\alpha]+t [\beta]) \in H_1(S)$ associated to $y=-s/t \in \hQ$. 
Such an element $\sigma_y$ is unique up to pre-composition with 
a power of the Dehn twist $\tau$ around $[\alpha]$, 
and we take and fix one of them. 
Then we obtain a map 
$$
 \B_\fty \sqcup \B_\fty^* \to  \B_y \sqcup \B_y^*
$$
by defining $[\rho] \mapsto \sigma_y \ten [\rho]$. 
Recall from Section 5.2 that the map 
$$
\Phi_\fty:\M \sqcup \M^* \to \B_\fty \sqcup \B_\fty^*
$$ 
is defined by $\Phi_\fty(\mu)=[\brho_\mu]$. 
Using these  maps, we define a map 
\begin{eqnarray*}
\Phi_y:\M \sqcup \M^* \to \B_y \sqcup \B_y^*
\end{eqnarray*}
by 
$\Phi_y(\mu):=\sigma_y \ten \Phi_\fty(\mu)$. 
Note that the map $\Phi_y$ is unique up to pre-composition with 
a power the translation $z \mapsto z+2$.  
We set 
$$
\B_y(p):=\Phi_y(\M(p)),  \quad    \B_y^*(p):=\Phi_y(\M^*(p))
$$
for integers $p \ge 1$. 

We denote by $\bd \D$ and $\inte(\D)$ 
the boundary and the interior of $\D$ in $\RR(S)$, respectively. 

\begin{defn}\label{bumpdef}
The space $\D$ is said to {\it self-bump} at 
$[\rho] \in \bd \D$ 
if there exists a neighborhood $U$ of $[\rho]$ 
such that for every neighborhood $V \subset U$ of $[\rho]$  
the intersection 
$V \cap \mathrm{int}(\D)$ is disconnected. 
\end{defn} 

We define the following two subsets of $\bd \D$: 
\begin{eqnarray*}
\bd^\bump \D&=&\{[\rho] \in \bd \D \,|\, \D \text{ self-bumps at } [\rho]\},  \\
\bd^\exotic \D&=&\{[\rho] \in \bd \D \,|\, 
\exists \, \text{an exotically convergent sequence with limit } [\rho]\}. 
\end{eqnarray*}

\begin{thm}\label{bump}
We have 
$$
\bd^\bump \D=\bd^\exotic \D=
\bigsqcup_{y \in \hQ}\left(\B_{y}(1) \sqcup \B_{y}^*(1)\right). 
$$
\end{thm}

We divide the proof of 
Theorem \ref{bump} into the following two Lemmas 
\ref{exolim} and \ref{bumpexo}; 
the former is a consequence of Theorems \ref{tan}  and \ref{limit}, 
and the latter is a consequence of 
Minsky's ending lamination theorem (Theorem \ref{Min}). 
In fact, we make use of Theorem \ref{Min} 
 to show that any self-bumping phenomena are induced by exotically convergent sequences. 
We remark that McMullen \cite{Mc1} used complex projective structures 
to show that any exotically convergent sequence 
induce self-bumping of $\D$ for general surface $S$. 

\begin{lem}\label{exolim}
$\bd^\exotic \D=\bigsqcup_{y \in \hQ}\left(\B_{y}(1) \sqcup \B_{y}^*(1)\right)$. 
\end{lem}

\begin{proof} 
Let $[\rho] \in \bd^\exotic \D$. 
Then by definition 
there is a convergent sequence $Q(x_n,y_n) \to [\rho]$ in $\D$
such that $(x_n,y_n)$ converges to some $(x_\fty, x_\fty) \in \Delta$. 
We have $x_\fty \in \hQ$ from Theorem \ref{Oh}. 
Our goal is to show that 
$[\rho] \in \B_{x_\fty}(1) \sqcup \B_{x_\fty}^*(1)$. 
By changing the generators of $\pi_1(S)$, 
we may assume that $x_\fty=\fty$. 
One can see from Theorems \ref{tan} and \ref{limit} that 
the limit $[\rho]$ of the exotically convergent sequence $Q(x_n,y_n)$ equals 
$[\brho_\xi]$ with $\xi=(p+1)\mu-p\ov\nu+2q$ 
for some $\mu, \nu \in \M$ and integers $p, q$ with $p \ne 0,-1$. 
It then follows 
from Lemma \ref{subset} (2)  that $\xi \in \M(1) \sqcup \M^*(1)$. 
Thus we obtain 
$[\rho]=[\brho_\xi] \in \B_\fty(1) \sqcup \B_\fty^*(1)$. 
 
To show the converse, 
it is enough to show that for a 
given $[\rho] \in \B_\fty(1)=\Phi_\fty(\M(1))$, 
there is an exotically convergent sequence with limit $[\rho]$. 
Take  $\xi \in \M(1)$ so that $[\rho]=[\brho_\xi]$. 
Then there exist $\mu,\nu \in \M$ such that $\xi=2\mu-\ov{\nu}$,  and we let 
$u, v$ be inverse images of $\mu,\nu$ via the map 
$\Psi:\bHH \sm \{ \fty \} \to \M$. 
Then the sequence $Q(\tau^n u,\tau^{2n} v)$ converges to $[\brho_\xi]$ as $n \to \fty$ 
by Theorem \ref{limit}. 
This completes the proof. 
\end{proof}

\begin{lem}\label{bumpexo}
$\bd^\bump \D=\bd^\exotic \D$. 
\end{lem}

\begin{proof}
We begin by fixing some notation of subsets of $\bHH \times \bHH$. 
We regard $\bHH \times \bHH$ as a subset of $\hC \times \hC$, and 
 define a distance 
between two points $(x,y), (x',y')$ in $\bHH \times \bHH$ by 
$$
d((x,y), (x',y'))=\max\{(d_s(x,x'), d_s(y,y')\}, 
$$
where $d_s(\cdot,\cdot)$ denotes the spherical distance in $\bHH \subset \hC$. 
(Although we do not make essential use of this metric, 
it is useful to understand precisely the arguments below.)  
Let  $[\rho_0]=Q(x_0,y_0) \in \bd \D$,  
which will be taken in $\bd^{\bump}\D$ or $\bd^\exotic \D$ in the sequel. 
Let $\ep_0=d((x_0,y_0),\Delta)=\inf_{z \in \bd \HH}d((x_0,y_0),(z,z))$ and take 
an open neighborhood ${\cal N}(\Delta)$ of $\Delta$ in $ \bHH \times \bHH$ as 
$$
{\cal N}(\Delta)=\{(x,y) \in \bHH \times \bHH  \,|\, d((x,y),\Delta) <\ep_0/2\}. 
$$
We set 
\begin{eqnarray*}
\ov{{\cal N}(\Delta)}&=&
\{(x,y) \in \bHH \times \bHH  \,|\, d((x,y),\Delta) \le \ep_0/2\}, \\
\bd {\cal N}(\Delta)&=&\ov{{\cal N}(\Delta)} \sm {\cal N}(\Delta)
=\{(x,y) \in \bHH \times \bHH  \,|\, d((x,y),\Delta) =\ep_0/2\}. 
\end{eqnarray*}

We first show that 
$\bd^\bump \D \subset \bd^\exotic \D$. 
Suppose that  $[\rho_0]=Q(x_0,y_0) \in \bd^\bump \D$, 
and let take an open neighborhood  ${\cal N}(\Delta)$ 
of $\Delta$ as above.  
We set 
$$
K=(\bHH \times \bHH) \sm {\cal N}(\Delta), 
$$
and let 
$$
\inte(K)=(\HH \times \HH) \sm \ov{{\cal N}(\Delta)} 
$$
be the interior of $K \subset \hC \times \hC$. 
Since $K$ is compact and since $\D$ is Hausdorff, 
it follows form Theorem \ref{Min} that 
the map 
$$
Q|_K: K \to \D
$$
is a homeomorphism onto its image. 
We remark  that 
$Q(\inte(K))=\inte(Q(K))$, 
which can be seen by checking that 
$\inte(K)$ is mapped by $Q$ into $\inte(Q(K))$ 
and $K \sm \inte(K)$ is mapped into $Q(K) \sm \inte(Q(K))$. 

It is  essential to observe that 
$Q(K)$ does not self-bump at $[\rho_0]$, 
whereas $\D$ does. 
In fact, first observe that 
$K$ does not self-bump at $(x_0,y_0)$; 
that is, for any neighborhood $U$ of $(x_0,y_0)$ in $K$ there is a 
neighborhood $V \subset U$ of $(x_0,y_0)$ such that 
$V \cap \inte(K)$ is connected.  
This can be seen precisely as follows: 
One can choose $\ep>0$ so that 
\begin{eqnarray*}
V=\{(x,y) \in \bHH \times \bHH \,|\, d((x,y), (x_0,y_0))<\ep\}
\end{eqnarray*}
is contained in $U$. 
In this case 
\begin{eqnarray*}
V \cap (\HH \times \HH)=
\{ x \in \HH \,|\, d_s(x,x_0)<\ep\} \times \{ y \in \HH \,|\, d_s(y,y_0)<\ep\}
\end{eqnarray*}
is connected. 
Since $Q(K)$ is homeomorphic to $K$ and since $Q(\inte(K))=\inte(Q(K))$, 
$Q(K)$ does not also self-bump at $[\rho_0]=Q(x_0,y_0)$; 
that is, 
for any neighborhood $U$ of $[\rho_0]$ in $Q(K)$ there is a 
neighborhood $V \subset U$ of $[\rho_0]$ such that 
$V \cap \inte(Q(K))$ is connected.  

Now let $U_0$ be a neighborhood of $[\rho_0]$ in $\D$ such that 
for every neighborhood $V \subset U_0$ of $[\rho_0]$, 
$V \cap \inte(\D)$ is disconnected.  
Since $Q(K)$ does not self-bump at $[\rho_0]$, 
one can choose  a decreasing sequence 
$$
U_0 \supset V_1 \supset  \cdots \supset V_n \supset  \cdots
$$
of neighborhoods of $[\rho_0]$ in $\D$ which satisfy the following:  
\begin{enumerate}
\item $\bigcap_{n=0}^\fty V_n =\{[\rho_0]\}$, and 
\item $V_n \cap \inte(Q(K))$ is connected for each $n$.
\end{enumerate}
Since $V_n \cap \inte(\D)$ is disconnected for each $n$,  
there is a sequence $\{(x_n,y_n)\}$ in $(\HH \times \HH) \sm \inte(K)$
 such that $Q(x_n,y_n)$ converges to $[\rho_0]$ as $n \to \fty$.  
Pass to a subsequence so that $(x_n,y_n)$ converges to some point 
$(x_\fty,y_\fty)$ in $(\bHH \times \bHH) \sm K$. 
If $(x_\fty,y_\fty) \not\in \Delta$ then 
$Q(x_n,y_n) \to Q(x_\fty,y_\fty)=[\rho_0]=Q(x_0,y_0)$, which contradicts 
the injectivity of the map $Q$. 
Thus $(x_\fty,y_\fty) \in \Delta$,  and thus  
$Q(x_n,y_n)$ is an  exotically convergent sequence.  
Therefore we obtain $[\rho_0] \in \bd^\exotic \D$. 

We next show that 
$\bd^\exotic \D \subset \bd^\bump \D$. 
Suppose that  $[\rho_0]=Q(x_0,y_0) \in \bd^\exotic \D$, 
and let take an open neighborhood  ${\cal N}(\Delta)$ 
of $\Delta$ as above.   
Observe that the complement $(\HH \times \HH) \sm \bd {\cal N}(\Delta)$ 
consists of two connected components $O_1,O_2$, where 
\begin{eqnarray*}
&&O_1=\{(x,y) \in \HH \times \HH  \,|\, d((x,y),\Delta) >\ep_0/2\}, \\
&&O_2=\{(x,y) \in \HH \times \HH  \,|\, d((x,y),\Delta) <\ep_0/2\}.
\end{eqnarray*}

We now claim that $Q(\bd {\cal N}(\Delta))$ is closed in $\D$. 
In fact, since  $\bd {\cal N}(\Delta)$ 
 is closed  in $\bHH \times \bHH$, it is compact. 
 Therefore the image $Q(\bd {\cal N}(\Delta))$ of 
 $\bd {\cal N}(\Delta)$ by the continuous map $Q$ 
 is also compact. 
Since $\D$ is Hausdorff,  we see that $Q(\bd {\cal N}(\Delta))$  is closed in $\D$.

Therefore we can choose  a neighborhood $U_0$ of $[\rho_0]$ 
so that 
$U_0 \cap Q(\bd {\cal N}(\Delta))=\emptyset$. 
Then, for any neighborhood $V \subset U_0$ of $[\rho_0]$, we have 
$$
V \cap \inte(\D)=(V \cap Q(O_1)) \cup (V \cap Q(O_2)). 
$$
Since $Q$ is a homeomorphism from $\HH \times \HH$ onto $\inte(\D)$, 
both $Q(O_1)$, $Q(O_2)$ are open,  and thus $V \cap Q(O_1)$ and 
$V \cap Q(O_2)$ are open. 
Since $[\rho_0]=Q(x_0,y_0) \in \bd^\exotic\D$, 
there is a sequence $\{(x_n,y_n)\}$ in ${\cal N}(\Delta)$ such that  
$(x_n,y_n)$ converges to a point in $\Delta$ and $Q(x_n,y_n)$ converges to $[\rho_0]$.  
By applying diagonal argument if necessary, we may assume that all $(x_n,y_n)$ lie in 
$O_2={\cal N}(\Delta) \cap (\HH \times \HH)$. 
Therefore we have $V \cap Q(O_2) \ne \emptyset$. 
On the other hand, it is obvious that  $V \cap Q(O_1) \ne \emptyset$.  
Thus $V \cap \inte(\D)$ is disconnected,  and thus 
we obtain $[\rho_0] \in \bd^\bump \D$. 
\end{proof}

\bigskip

\begin{flushleft}
Graduate School of Mathematics, \\
Nagoya University, \\
Nagoya 464-8602, Japan  \\
\texttt{itoken@math.nagoya-u.ac.jp}  
\end{flushleft}


\begin{thebibliography}{Brock} 

\bibitem[AC]{AC}J.~W.~Anderson and  R.~D.~Canary.  
{\it Algebraic limits of Kleinian groups which rearrange the pages of a book.}  
 Invent. Math. {\bf 126} (1996),  no. 2, 205--214. 

\bibitem[Bo]{Bo}F.~Bonahon. 
{\it Bouts des vari\'et\'es hyperboliques de dimension $3$.} 
Ann. of Math. (2) {\bf 124} (1986),  no. 1, 71--158. 

\bibitem[Bow]{Bow}B.~H.~Bowditch.
{\it Markoff triples and guasifuchsian groups.}
Proc. London Math. Soc.  {\bf 77} (1998), 697--736. 

\bibitem[Brock]{Brock}J.~F.~Brock. 
{\it Iteration of mapping classes on a Bers slice: 
examples of algebraic and geometric limits of hyperbolic $3$-manifolds.}   
Lipa's legacy (New York, 1995), 81--106, Contemp. Math., {\bf 211}, 
Amer. Math. Soc., Providence, RI, 1997. 

\bibitem[Brom]{Brom}K.~W.~Bromberg. 
{\it The space of Kleinian punctured torus groups is not locally connected.} 
to apper, Duke Math. J. 

\bibitem[BH]{BH}K.~W.~Bromberg and J.~Holt.  
{\it Self-bumping of deformation spaces of hyperbolic 3-manifolds.} 
J. Differential Geom. {\bf 57} (2001), no. 1, 47--65. 

\bibitem[Ca1]{Ca1}R.~D.~Canary. 
{\it Pushing the boundary.}  
In the tradition of Ahlfors and Bers, III, 109--121, 
Contemp. Math., {\bf 355}, Amer. Math. Soc., Providence, RI, 2004. 

\bibitem[Ca2]{Ca2}R.~D.~Canary. 
 {\it Introductory Bumponomics: the topology of deformation spaces of hyperbolic 3-manifolds}, 
 in Teichm\"{u}ller Theory and Moduli Problem, ed. by I. Biswas, R. Kulkarni and S. Mitra, Ramanujan Mathematical Society, 2010, 131-150. 

\bibitem[EM]{EM}D.~B.~A.~Epstein and A.~Marden. 
{\it Convex hulls in hyperbolic space, a theorem of Sullivan, and measured pleated surfaces.}  
Fundamentals of hyperbolic geometry: selected expositions, 117--266, 
London Math. Soc. Lecture Note Ser., {\bf 328}, Cambridge Univ. Press, Cambridge, 2006. 

\bibitem[J{\o}]{Jo}T.~J{\o}rgensen. 
{\it On discrete groups of Mobius transformations.} 
 Amer. J. Math. {\bf 98} (1976), no. 3, 739--749. 

\bibitem[JM]{JM}T.~J{\o}rgensen and A.~Marden. 
{\it Algebraic and geometric convergence of Kleinian groups.}  
Math. Scand. {\bf 66} (1990), no. 1, 47--72. 

\bibitem[KS]{KS}
L.~Keen  and C.~Series.  
{\it Pleating coordinates for the Maskit embedding of the Teichm\"{u}ller space of punctured tori.} 
Topology {\bf 32} (1993), no. 4, 719--749. 

\bibitem[KT]{KT}S.~P.~Kerckhoff and W.~P.~Thurston.  
{\it Noncontinuity of the action of the modular group at Bers' boundary of Teichmuller space.}   
Invent. Math. {\bf 100} (1990), no. 1, 25--47. 

\bibitem[Ma]{Ma}
A.~Marden. 
Outer circles An introduction to hyperbolic 3-manifolds. 
Cambridge University Press, 2007.  

\bibitem[Mc1]{Mc1}C.~T.~McMullen. 
{\it Complex earthquakes and Teichmuller theory.}  
 J. Amer. Math. Soc. {\bf 11} (1998), no. 2, 283--320. 

\bibitem[Mc2]{Mc2}C.~T.~McMullen. 
{\it Hausdorff dimension and conformal dynamics. I. Strong convergence of Kleinian groups.}  
 J. Differential Geom. {\bf 51} (1999), no. 3, 471--515. 

\bibitem[Mi]{Mi}Y.~N.~Minsky.  
{\it The classification of punctured-torus groups.}  
Ann. of Math. (2) {\bf 149} (1999), no. 2, 559--626. 

\bibitem[Oh1]{Oh1}K.~Ohshika. 
{\it Divergent sequences of Kleinian groups.} 
The Epstein birthday schrift, 419--450,  
Geom. Topol. Monogr., {\bf 1}, Geom. Topol. Publ., Coventry, 1998.

\bibitem[Oh2]{Oh2}K.~Ohshika. 
{\it Divergence, exotic convergence and self-bumping in quasi-Fuchsian spaces.} 
Preprint, arXiv:1010.0070. 

\end{thebibliography}
\end{document}